\documentclass[10pt]{article}
\usepackage[T1]{fontenc}
\usepackage{amssymb}
\include{amssymb}
\usepackage{amsmath}
\usepackage{graphicx}

\newcommand{\Rn}{\mathbf{R}^{n}}

\begin{document}
\title{Some remarks on Schr\"odinger operators}
\author{T. Dahn (Lund University)}

\maketitle 

\section{Introduction}
Starting with a mixed model $f(\zeta) \rightarrow F(\gamma)(\zeta) \rightarrow \gamma \rightarrow \tilde{\gamma}$
(\cite{Dahn17}), we have the following problem. Given a potential $\widehat{V}$ such that $\log \widehat{V} \in L^{1}$,
there is a corresponding movement of the base for the symbol ideal, such that $F(U \gamma)=F(\gamma) + \widehat{V}$
and conversely, given a movement $U$ we can determine $\widehat{V}$. In particular we consider the spiral
movement $U_{S}$ and the corresponding potential.
The Schr\"odinger opertors in this article, give a global model that is not a normal model.

\subsection{invariance principle}
The invariance principle (\cite{Dahn17}) is relative a continuous mapping $J : \gamma_{1} \rightarrow \gamma_{2}$ with
$\gamma_{1}=\gamma_{2}$ on a set $\Delta$, where $\gamma_{j}$ analytic and $\gamma_{2}$ reduced with
respect to $\gamma_{1}$. We assume $F (J \gamma_{1})=F(e^{v} \gamma_{2})$ for $v \in L^{1}$ and the 
boundary is assumed very regular (\cite{Dahn13}).
The invariance principle is given $\gamma_{1}$ closed, $U$ a movement of $\gamma_{1}$ then $V=JU$,
where $V$ is a movement in $\gamma_{2}$, conversely $U=J^{-1} V$.

Consider the mixed model, $F(\gamma) \rightarrow \gamma_{1} \rightarrow \gamma_{2} \rightarrow \zeta$
where $\gamma_{2} \in (I_{1}^{\bot_{E}})$ and $J U = V J$, for a movement $U$ with generator $R$.
If between two points $p,q$, the generator $R$ has constant angle, $\frac{d y'}{d x'}=\frac{d y}{d x}$, we
consider the continuation $\tilde{R}$ to the infinity. The movement changes character when
$d y/d x-1$ changes sign ($ \leq 0,\geq 0, =0$). 
Assume $y=y(x)$ and $\eta(x)=y(x)/x$ and $J \gamma_{1}=\gamma_{2}$. The condition $\frac{d y_{2}}{d x_{2}}=\frac{d y_{1}}{d x_{1}}$
means that a decomposition in a composite movement, can be given independent of $J$.
Given a composite movement, $J(x,y) \rightarrow (\tilde{x},\tilde{y})$ and 
$A(\tilde{x},\tilde{y})$ a linear movement, then we have existence of $U(x,y)$, such that
$J U=A$. If the movement is not dependent on $P_{j}=(x,y)$, $j=1,2$, we have $A=V$.

The de-singularization with $U$ is not present for $V$, that is we may have not-discrete intersection
for invariant axes. We consider instead the problem $F(\gamma) \rightarrow \gamma \rightarrow \gamma' \rightarrow \zeta$
where $\gamma'$ is polynomial. (\cite{Nishino75}) Note that $\int_{\Omega} \gamma' d x=0$
implies $m \Omega=0$ (Hurwitz theorem). When we consider $JU=V$ where $U=V$ over the set of lineality, 
we note that $\Omega_{1}=\{ \zeta \quad U \gamma - \gamma \}$
and $\Omega_{2}=\{ \zeta \quad V \gamma - \gamma \}$, have different properties under iteration.
$\Omega_{1}$ is invariant and $\Omega_{2}$ decreasing for iteration.  Further, note
when $F$ is a localizer, the kernel to respective $F,\tilde{F}$ have different properties.

Assume $F(U \gamma)=F(\gamma) + \widehat{V}$, such that $\widehat{V}=0$ is symmetric in a neighbourhood
of a point with respect to $\frac{y}{x}=\eta \rightarrow 1/\eta$. When $U$ is surjective (projective) and $F$ has trivial kernel, 
$\widehat{V}=0$ is given by $U^{\bot_{E}} \gamma$. When $\eta$ in a semi-analytic set, where we have change of base,
we are considering $\{ \mu \leq \eta \leq \lambda \}$, for constants $\mu,\lambda$. On points $C$, 
where we have a change of base $\eta \rightarrow 1/\eta$,
we have $\frac{d \eta}{d x}=0$ implies $\eta=const$, that is the change $\eta \rightarrow 1/\eta$
is a jump.

The invariance principle, means that we can uniquely determine the movement by reverting
$U \eta \rightarrow V \eta \rightarrow V \eta^{N}$.
The model so far assumes
$J$ preserves a constant sign, that is $0 \notin \sigma(J)$. 

\newtheorem{prop2}{Proposition}[section]
\begin{prop2}
 The representation $F(U^{\bot_{E}} \gamma)=\widehat{V}(\zeta_{t})$ defines $\widehat{V}$ as dependent on the movement
$U$ and on $F$ and as independent on $\gamma$, the base. 
\end{prop2}

When $F(U^{\bot_{E}} \gamma)=\widehat{V}$, we have that $\widehat{V}=0$ if either $U^{\bot_{E}} \gamma$ is 
reflection axes or in the kernel to $F$. For the spiral all points are in $C$ and simultaneously the 
movement is monotonous, thus it must be situated in the kernel to $F$.

When the set for change of character $C$ is a discrete set
in the finite plane, we are
considering consecutive simple movements. Note that change of axes does not necessarily imply change of character.
When $x \rightarrow \eta(x)$ bijective, we can consider regularly moving reflection axes, in this case 
the change of character is through $\eta=1$. Otherwise, through scaling of the hyperboloid, we can assume
change of character of movement is arbitrarily close to $\eta=1$.

When for $v=1/\eta$, the set $C \cap \{ \eta + i v \}$
are points, we have simple movements, when it is intervals, we have spiral like movements.
We assume, where the movement changes character, we have $\widehat{V}(\zeta_{t}) \rightarrow \zeta_{t}$
continuous, where $t$ is parameter for the movement. Further, when $d U_{1} = \alpha d U_{2}$, we have $\alpha \frac{d y}{d U_{1}} = \frac{d y}{ d U_{2}}$,
where $\alpha$ regular, and we assume the mapping $\gamma \rightarrow U$ is differentiable (modulo monotropy).

\subsection{Representation using parametrises}
Assume $f(D) + V(x)$ is defined such that $(f(D) + V)^{\wedge}=f(\zeta) + \delta_{0}$. The properties
of $V$ are determined by the topology. Assume that $E$ is a parametrix to the operator $f(D)$
so that $E f(D)=\delta_{0} - \gamma$, for $\gamma \in C^{\infty}$. We then have that if $V(x) \in \mbox{ ker }E$,
then $E (f(D) + V(x))=E f(D)=\delta_{0} - \gamma$. Given that $E_{N} = \{ 0 \}$ (iteration)
if we let $V_{N}$ be the potential corresponding to $f^{N}(D)$, then obviously $V_{N}$ does not
have support in $0$. Thus, we have $V_{N} \neq const.$. If the kernel is taken modulo $C^{\infty}$,
we can assume $V_{N} \in C^{\infty}$. 

Assume $E_{0}$ such that  $PE_{0}=\delta$ and $E$ such that $PE=\delta_{0} - \gamma$,
for $\gamma \in C^{\infty}$ where $E \sim E_{0} + R$. We assume $PE=PE_{0} + PR$, where $PR \in C^{\infty}$
(parametrix method). Note that $E : \mathcal{D}' \rightarrow \mathcal{D}^{' F}$ and $E_{0} : \mathcal{D}' \rightarrow \mathcal{D}'$.
Assume $E$ parametrix to $f(D)$, such that $(E-I)f \in C^{\infty}$ outside the kernel to $E$ and we have existence
of $P(D)$ such that $P(D)E(f)=f + C^{\infty}$. 
If we assume $P(D) E=E P(D)$, that is $E$ is two-sided, this means that $\big[ f,I \big]=\big[ I,f \big]$ modulo $C^{\infty}$,
outside the kernel to $E$. Further $PEPE(f) \sim P^{2}E_{2}(f)$, why we have that $f$ is algebraic on the domain to $E$. 
Further, assume $F_{1}$ corresponds to a parametrix to the base $\gamma$ and $F(\gamma)=\big[ f, F_{1} \big](\gamma)$, where $F$
is localizer to $f$. That is $F_{1} \gamma = I$ (modulo $\widehat{C}^{\infty}$, that is corresponding to regularizing action) and $\big[ f,I \big]=F(\gamma)$.

The movement can be characterized through $<\widehat{V}, W \phi>=<U^{\bot_{E}} \gamma, {}^{t} F(W \phi)>=0$
where we assume $<\phi,\gamma> =1$, for $\phi \in H_{m}$ and $\gamma \in (I)$ an ideal for the base.
When $f \in L^{1}$, we have isolated singularities and we can use monotropy to chose $f_{1}$ holomorphic
close to the singularity and $f$ continuous such that $f \sim_{m} f_{1}$ (\cite{Dahn13}).
The condition $\log \eta \in L^{1}$ means particularly for $\eta = e^{\phi}$, that $\phi \sim_{m} \phi_{1}$
where $\phi_{1}$ is locally analytic, thus $\eta = const \sim$ an analytic set.

 \subsection{On presence of trace}
Presence of a trace implies presence of lineality $\Delta$. Further, translation invariance
for $g_{\zeta}(x)=\int E_{\theta(\zeta)}(x,y) e^{i <y,\zeta >} d y$ implies a trace for $E_{\theta}$.
Starting with a contact transform $y \rightarrow \xi$, this can be continued
in one variable, to $y \rightarrow \zeta$, as in the one-dimensional Laplace transform. Thus when
$g_{\zeta}$ holomorphic in $\zeta$ and
$g_{\zeta} \rightarrow 0$ on a ray in $\zeta$, we have 
$g_{\zeta}=0$ on a sector in $\zeta$. 
Given that translation is algebraic over $g$, the sector is a full disc.
When $E_{\theta}$ parametrix to a differential operator, this means presence of kernel,
that is $E_{\theta(\zeta)}(\varphi)=0$ for $\varphi \in \mathcal{D}(Y)$, where $Y$ is non-trivial.
.  

When $ {}^{t} U I_{E}(\phi)=\int E({}^{t} U_{1} x,y) \phi(y) d y$ and $I_{E}(U \phi)=\int E(x,y) U \phi(y) d y$,
let $d y'=\rho d y$, we have $I_{E}(\phi)=\int E(x,{}^{t} U_{1}^{-1} y') \phi(y') d y'/\rho$. When $I_{E}$
is not nuclear, for instance $U^{\bot_{E} \bot_{E}} \neq U$ or ${}^{t} U_{1}^{\bot_{E}}$ not closed, when $U^{\bot_{E}}$ closed,
we do not have ${}^{t} U I_{E} = I_{E} U$. When $I_{E}$ is nuclear, we have $<{}^{t} U^{\bot_{E}} I_{E}(\phi),\psi>=
<E, U^{\bot_{E}} \psi \otimes \phi>$ and $<I_{E}(U^{\bot_{E}} \phi),\psi>=<E,\psi \otimes U^{\bot_{E}} \phi>$, thus
equivalence requires a symmetric domain for $I_{E}$, in particular the kernel must be symmetric.

Assume $H_{V}=-\Delta + V$ corresponds to a movement on $H_{0}$. Thus $H_{V}=H_{0}$ corresponds to invariant points.
Assume $E$ a two-sided parametrix, such that $H_{V} E = E H_{0}$ implies $V E =0$.
When $V$ corresponds to a movement, invariant points are in $\mbox{ ker }E$.
Note that $E$ can be selected as hypoelliptic over the kernel, when
$H_{V}$ is partially hypoelliptic. 
(cf. \cite{Melin87})

\subsection{On a transmission property}
Note that when a mapping bijective, maps zero-lines on to zero-lines, it is projective.
Concerning $\sigma(J) \neq \emptyset$, we must have a set $\Delta$, where $\gamma_{1} \sim \gamma_{2}$,
in particular $\Delta(\gamma_{1})=\Delta(\gamma_{2})$. Note that $\gamma_{j}$, $j=1,2$ have
different properties under iteration. Given $J : 0 \rightarrow 0$, we have over $\Delta$,
$\tau J = J \tau$ implies $J (\tau - 1) \simeq (\tau - 1)$, where $\tau$ is translation in domain for $\gamma_{j}$. 

\newtheorem{def3}[prop2]{Definition}
\begin{def3}
Define $\Xi(x,y)=\big[ F(x,\eta), {}^{t}F(y,v) \big]$, where $\eta=y/x$ and $v=x/y$. When $\eta \rightarrow v$
is projective, $\Xi(x,y) \sim F(x,y) \int d \mu(\eta,v)$ as a Schwartz kernel.
\end{def3}

Assume $F(x,y) = F_{1}(x,\eta)$ and ${}^{t}F(x,y)=F_{2}(y,v)$. Given that $(\eta,v) \rightarrow (v,\eta)$
is projective and further $\eta^{2}$ reduced, then $\int_{\Omega} d \eta=0$ implies $\int_{\Omega} \eta^{2} d v=0$
and $v (\Omega)=0$ (measure zero). Further, $\frac{d F}{d x} \sim \frac{d}{d x} \big[ F_{1},{}^{t} F_{2} \big]$
and $\frac{d F}{d y} \sim \frac{d}{d y}\big[ F_{1},{}^{t} F_{2} \big]$. Geometrically the representations of $F$
are equivalent.
We extend the integration domain as above to $(x,\eta;y,v)$ and consider $(F(x,y),M(\eta,v))$ where
$M(\eta,v)$ is integrand corresponding to $d \mu$.

 The transmission property can be discussed relative
movements, if it is present for all movements, the property is ''global``.
When the movement is mixed, we assume $\log U_{1} U_{2} \eta \in L^{1}$, where we have isolated singularities.
Consider $ (\eta,v) \in (I) \times (I)^{\bot_{E}}$.
When $U_{1} : (I) \rightarrow (I)$ and ${}^{t} U_{2} v \in (I)^{\bot_{E}}$ iff $v \in (U I)^{\bot_{E}}$.  We have
${}^{t} U_{2} : (U I)^{\bot_{E}} \rightarrow (I)^{\bot_{E}}$. Note that when we have projectivity, $(I^{\bot_{E} \bot_{E}}) \cap (I)= \{ 0 \}$
Note that $\frac{d}{d x}\big[ U_{1},U_{2} \big] = \frac{d U_{1}}{d U_{2}} \frac{d U_{2}}{d x} = \alpha \frac{d U_{2}}{d x}$
where $\alpha \neq 0$ is regular. 
Note that the domain where $U_{1} U_{2} = U_{2} U_{1}$ is symmetric has an interpolation property.
The problem that ${}^{t} M \sim M$ implies monodromy in both ends, is dealt with by assuming symmetry
in $H'$, but not in $H$ (or $L^{2}$)

Points where the movement changes character $C$ are such that ${}^{t} M(\eta,v)= \lambda M(\eta,v)$
and $M(U^{\bot_{E}} \eta,v)(\zeta)=\widehat{V}(\zeta)$, that is we may adjust the axes with $\lambda$,
without changing character.
Note that $C$ has a connected set of $(x,y)$ for a spiral and that $\widehat{V} = const$ over $C$.
 Given 
$\mid U_{1} \gamma \mid \sim \mid U_{2}^{\bot_{E}} \gamma \mid$, we assume these points of $C$ are on $H_{m}$ Note that the spiral can not 
completely be represented on $H_{m}$. 

When $C$ is discrete, the movement can be completely determined by the hyperboloid, that is it is assumed independent of $J$.
When $\frac{d y}{d x}$ has isolated zero's, the same holds for $\frac{d y}{d x}-1$. 
Note that $\frac{d y}{d x} - 1=\frac{d}{d x}(y-x)$ and when $\frac{d}{d x}e^{y-x}=1$, $\frac{d y}{d x} < 1$
when $x < y$.
When $(y-x) \in L^{1}$ and $x$ reduced in $\zeta$, we have $x(1-\eta) \in L^{1}$ and using Nullstellensatz,
$(1 - \eta)=\eta^{\bot_{E}} \in L^{1}$. 
Note that $\frac{d}{d x}e^{y-x} = (\frac{d y}{d x})^{\bot_{E}} e^{y-x}$. 

\newtheorem{lem1}[prop2]{Lemma}
\begin{lem1}
 A point where $u$ is both rotation and translation,
is written $v \sim \sqrt{u_{1}u_{2}}$ and $\frac{d v}{d t} > 0$, that is monotonous. Assume $v^{2} \sim u_{1}u_{2}$
and $(d w)^{2} \sim d v^{2}$.  The movement looked for 
is $\int d w \sim W$ and $u_{1}u_{2} \leq W \leq u_{2} u_{1}$. We assume $W$ symmetric relative $\eta,1/\eta$
Note that $(\frac{d}{d t} w)^{2} \sim \frac{d}{d t} u_{1}u_{2} \sim \frac{d}{d t} u_{2} u_{1}$. When $W$ is a.c
$d w=0$ implies $\int d w=const$ that is $W$ is completely symmetric.
\end{lem1}

Note that $v^{2}$ a.c. (absolute continuous) does not imply $v$ a.c.
Assume $d v^{2}$ a positive measure, defining a bounded operator on a Hilbert-space, 
then there is a positive measure $d w$ such that $d v^{2} \sim (d w)^{2}$.
If $d C$ is a measure, defining a bounded operator on a Hilbert-space, such that $ \big[ dC,d v^{2} \big] \sim \big[ d v^{2}, d C \big]$, then
$\big[ d C, d w \big] \sim \big[ d w , d C \big]$.
Consider $d w \overline{d w} \sim d w( z, \overline{z})$ or more generally
$d w (\eta, v)$. We assume that $(\eta,v)$ has compact sub-level surfaces, why the measure can be defined
on compact sets. When $d w \equiv 0$ over $(\eta,v)$, we can assume $W$ (a.c.) completely symmetric over 
reflection axes represented in $(\eta,v)$.
Note that when $(d w)^{2} \equiv 0$ and a.c., we have $v^{2} \sim const$, but when for instance $d w = 0$ a.e.
we do not have $w=const$ (\cite{Riesz56}, Ch. 2, Section 24)

A pure differential is such that $(w + i w^{*})^{*}=-i (w + i w^{*})$.
When $\varphi^{*}=-i \varphi$, we have that $\varphi$ is pure. 
When $\varphi$ is closed and $\varphi^{*}=-i \varphi$, then $\varphi$ is analytic.(\cite{AhlforsSario60}) 

Concerning $\big[ U,I \big]=\big[ I,U \big]$, we note that when we have unbounded sub-level sets
for the symbol, the corresponding spectral kernel is a distribution. That is, when $F$ is not corresponding to 
very regular action,
for instance $\mbox{ ker }F$ non-trivial, algebraicity is not established.  
When $U \gamma \in \mathcal{D}_{L^{1}}$, we can define $T F=FU$ in $\mathcal{D}_{L^{1}}'$.
When $F(\gamma)=e^{\phi}$ and $T F(\gamma)=e^{t \phi}$ and $t$ algebraic, we assume
$t \phi \in L^{1}$.

\section{The kernel theorem}
\subsection{Schwartz kernel theorem}
Assume $F$ has Schwartz kernel and is algebraic in the sense that $\big[ F, U \big] = \big[ U,F \big]$, 
then $<F,{}^{t} U \varphi \times \phi>=<F,\varphi \times U \phi>$.
When $F$ is symmetric we have $U \varphi \times \phi \simeq \varphi \times U \phi$. Thus, when $F$ is algebraic over $U$
and symmetric, we have $U \rightarrow {}^{t} U$ is bijective over $\mathcal{D}$. 

When $\widehat{V}=F(U^{\bot_{E}} \gamma)$ and $\widehat{V}=0$ implies $\mid \eta \mid \leq 1$, we have
semi-bounded support for $F U^{\bot}$. When $\widehat{V} \sim {}^{t} \widehat{V}$, that is $\mid 1/ \eta \mid \leq 1$,
we have bounded support. Consequently the set $\{ \zeta \quad \lambda \leq \mid \eta \mid \leq \mu \}$ is unbounded in  both
ends. A sufficient condition for a regular spectral kernel is that $\eta$ corresponds to a hypoelliptic operator, in this case the 
set above is compact (\cite{Nilsson72}).
Note that $\eta(x)$ is algebraic (in $x$) iff $y(x)$ algebraic (in $x$).
Further, the condition $\mid x/y \mid \rightarrow 0$ as $\mid \zeta \mid \rightarrow \infty$, implies that
$\{ \zeta \mid \eta \mid \leq \lambda \}$ bounded (and closed). 

Using Lie's method of integration we can write
$F(x,y) = G(x , \eta(x))$ with
$\lambda y(\lambda x)/\lambda x=y(\lambda x)/x$ and if $y(\lambda x)=\lambda x$
we have $\lambda \eta(\lambda x)=\lambda \eta(x)$. Thus we consider $(x, y) \rightarrow (x,\eta) 
\rightarrow \lambda (x',\eta') \rightarrow  \lambda (x',y')$,
where the mapping $J $ maps the hyperbolic space $(x,y)$ on to euclidean space $(x',y')$.
Thus $(x,y) \rightarrow  \lambda (x',\eta')$ corresponds to scaling. The set $\Sigma=\{ (x',y') \quad \mid \eta' \mid \leq \lambda \}$
has a corresponding set $\Omega=\{ \zeta \quad (x',y') \in \Sigma \}$. 
Note when $\zeta \rightarrow U \gamma(\zeta)$ is considered on $\gamma \in H_{m}$, it is assumed analytic.
In general it is only assumed continuous on $U^{\bot_{E}} \gamma$, when $\gamma \in (I_{1})$. 

A proper mapping $J$, is a continuous mapping, that has a continuous continuation to $(I_{1})^{*} \rightarrow (I_{2})^{*}$
(one-point compactification), such that $J (\infty_{1})=\infty_{2}$. 
The condition on ``collar points'' must be assumed in $(I_{2})$ as well as $(I_{1})$. The condition we use 
is $d (U_{1} - U_{2})=(U_{1} - U_{2})=0$ implies $\gamma=\gamma_{0}$, a point (cf. \cite{Lie96}). Any movement on $H_{m}$ has correspondent
sub-level sets in $\zeta$.  Thus, assume collar points 
in both $U$ and $JU=V$, then $\mid \eta'' \mid \leq \lambda$ has a correspondent $\mid U_{3} \eta' \mid \leq 1$
where $U_{3}$ is parabolic. That is a movement relative the $U_{3}$ axes, corresponds to scaling with $\lambda$.
Note that when the simple movements are taken in sequence, we have 
$\tilde{U}^{2} \sim (\tilde{U}_{1}^{\bot_{E}})^{2}=(U_{1})^{\bot_{E}}(U_{1}^{*})^{\bot_{E}}$ 
has a correspondent $\tilde{V}^{2} \sim (\tilde{V}_{1}^{\bot_{E}})^{2} =(V_{1})^{\bot_{E}} (V_{1}^{*})^{\bot_{E}}$,
where $U^{*}$ is reflection in a cylindrical domain.

The spiral is approximated by a sequential change of axes $\eta_{0} \rightarrow U_{1} \eta_{0}=\eta_{1}$
and $\eta_{1} \rightarrow U_{2} \eta_{1}=v_{1}$ and so on. When the displacement is made infinitely small,
we get the spiral movement, a connected symmetry set for $M$, which in this case is not analytic. 
 
 \subsection{Symmetry}
Consider $\{ (\eta,v) \quad P(\eta,v) \leq \lambda \}$ where $v \sim 1/\eta$.
Tarski-Seidenberg gives that $\{ \eta \quad P(\eta,v) \leq \lambda \}$ is semi-algebraic. 
When $P(\eta,v) \sim P_{1}(\eta) P_{2}(v)$ and $P_{2} \sim 1/Q(\eta)$,
$\{ P(\eta,v) \leq \lambda \}$ iff $\{ P_{1}(\eta)/Q(\eta) \leq \lambda \}$. Thus, the set is semi-algebraic
when $P_{1}/Q$ polynomial. The corresponding proposition for compact sub-level sets, requires $P_{1},P_{2}$
reduced. When ${}^{t} P \sim P$, $P$ is not necessarily algebraic. 
When $M(\tilde{U}_{1} \eta,v) \sim M(\eta, \tilde{U}_{2} v)$ and when ${}^{t} M  \sim M$, we could say that 
$M$ is symmetric with respect to $\tilde{U}_{1},\tilde{U}_{2}$ close to points in $C$. 

Note that if $\{ F(\gamma) < \lambda \}$ is semi-algebraic, we have that $\{ x \quad F < \lambda \} \cup \{ y \quad F < \lambda \}$
is semi-algebraic. Using the involution $y \rightarrow \frac{1}{y}$, when we have existence
of $Q$ polynomial such that $F(x,y) \sim Q(x,\frac{1}{y})$, when $y \rightarrow \infty$
(preserves constant value, \cite{Cousin95}), then $\{ \frac{x}{y} \}$ is semi-algebraic.

Assume $<\big[ \tau F, I \big] \phi,\psi>=<\big[ F, {}^{t} \tau I \big] \phi,\psi>$. If $\big[ F,I \big]=\big[ I , F \big]$,
we have that $<\big[ \tau I, F \big] \phi,\psi>=<\big[ F, {}^{t} \tau I \big] \phi,\psi>$ and if
$<\big[ \tau I, F \big] \phi,\psi>=<\big[ F, \tau I \big] \phi,\psi>$, we must have that $\tau \rightarrow {}^{t} \tau$
is bijective. Assume $\big[ U I , F \big] = \big[ F , U I \big]$, then 
$< F, {}^{t} U \varphi \times \phi > = < F, \varphi \times U \phi>$, thus when $U \rightarrow {}^{t} U$
bijective over $\mathcal{D}$, we have $F$ algebraic (and symmetric) over $U$.

Assume $\Phi(\frac{d y}{d x})=\eta$, the blow-up mapping.
The condition for determined tangent $\frac{d y}{d x}=\frac{y}{x}$ (\cite{Bendixsson01}), that is $\Phi \equiv 1$,
applied to $(U_{1} - U_{2})$ means through collar point condition, that singularities with determined tangent are points. Assume $\psi(\eta)=1/\eta$
and consider $\Phi=\lambda \psi$, which is usually a jump discontinuity. The mapping $\frac{d y}{ d x} \rightarrow \lambda \frac{d x}{d y}$
is dependent on $y=y(x) \rightarrow  x=x(y)$, which is dependent on involution. Over constant surfaces,
we have $(\frac{d y}{d x})^{2} \sim \lambda$. Note that $\eta^{2}$ polynomial does not imply $\eta$ polynomial.

\subsection{The sub-level sets}
Consider $F(U \gamma)-F(\gamma)=\widehat{V}$ and $\gamma=(x,y) \rightarrow \frac{y}{x}=\eta(x)$ where
$y=y(x)$ in a rotation surface. 
Define $(F(\gamma),M(\eta, 1/\eta))$, so that ${}^{t} M(\eta, 1/\eta) \sim M(1/\eta,\eta)$.
In the vicinity of a point where the movement changes character,
we assume ${}^{t}M \sim M$ and that on $C$, we have $\widehat{V}=0$. 

 Assume $F(U^{\bot_{E}} \gamma)(\zeta) \sim F(\gamma)(\zeta_{t})$,
where $U^{\bot_{E}} \rightarrow \zeta_{t}$ is continuous on $C$, the set for change of character, 
set of singularities for $\widehat{V}$.
When this set in $\zeta$ is not discrete, we have to assume regular approximations
in $\zeta_{t}$ corresponding to simple consecutive movements.

Consider $\{ \lambda < \mid F \mid < \mu \}$, when the set is compact in $(x,y)$, it is not necessarily
compact in $\zeta$. 
The condition $1/f \rightarrow 0$ in $\infty$, implies
$\mid \zeta \mid^{c} \leq C \mid f(\zeta) \mid$, outside a compact set for $\zeta$ in the real space.
Thus $\Omega = \{ f < \lambda \}$ for a constant $\lambda$ is compact. In the complex space we consider $f$ in $\mathbf{R}^{2 n}$.
When $\Sigma = \{ x \quad \eta(x) \leq \lambda \}$ and assuming $x$ continuous and reduced in $\zeta$, then $\Omega$
is compact iff $\Sigma$ compact.

Assume $\frac{1}{\eta + i \frac{1}{\eta}} \rightarrow 0$, as $\mid \zeta \mid \rightarrow \infty$
When $\mid M(\eta,1/\eta) \mid \leq C (\mid \eta + i 1/\eta \mid)$ and $f(\zeta)$ hypoelliptic, we have that $\eta + i 1/\eta$
reduced. When $\eta \rightarrow c$ as $\mid \zeta \mid \rightarrow \infty$, we can consider $\eta'=\eta-c$,
since this is only a scaling of the hyperboloid. 
When
$M(\eta,1/\eta) \sim  M(1/\eta,\eta)$ as $\mid \zeta \mid \rightarrow \infty$, the sets
$\lambda_{1} < \mid \eta \mid < \lambda_{2}$  may be unbounded in both ends. Note that Fredholm index is not constant
on these sets. 
When $\phi \in \mbox{ ker } {}^{t} M$, we have $\eta + i \frac{1}{\eta} \bot_{E} {}^{t} M \phi$, for any $\eta$,
that is when we have a non-trivial kernel for ${}^{t} M$, there is space for a spiral symmetric set.
Note the example when $\phi (\frac{1}{\eta}) \sim \frac{1}{ \phi}(\eta)$, then $\int M(\eta,1/\eta) \phi(\frac{1}{\eta}) \sim 
\int_{\mid \eta \mid > 1} M(1/\eta, \eta) \frac{1}{\phi}(\eta)$. Thus when $M$ has support
in a neighbourhood of $\infty$ (as a distribution), we must have $M \neq {}^{t} M$ .

However $\eta + i v$, $v \sim 1/\eta$
has compact sub-level sets and this domain is symmetric in the sense that $v + i \eta$ has compact sub-level sets
and over this set the movement can change character. 
Note that using the blow-up mapping $d y \rightarrow \eta d x$, $\eta dy \rightarrow \eta^{2} d x$ and
$d \eta =- \eta^{2} d v$. A transmission property in this context, is that the spiral can be rectified
to a line in finitely many steps.
According to the generalized moment problem (\cite{Riesz56}), 
$E_{0} \bot_{E} d \mu$ implies $C^{0} \bot_{E} d \mu$, 
assuming $d \mu$ of bounded variation and $E_{0} \subset C^{0}$,
given a separation property, for instance the Schwartz kernel to $F$ corresponding to very regular action in euclidean space. 

\subsection{Composition of functionals}
Consider a composition of movements as a functional $\mathcal{D}'(X \times Y)  \ni \tau_{\lambda \mu} \simeq \tau_{\lambda} \tau_{\mu} \in \mathcal{D}'(X) \times \mathcal{D}'(Y)$
given that $\tau_{\lambda} \tau_{\mu}=\tau_{\mu} \tau_{\lambda}$, we have an interpolation property. In general we have $\mathcal{D}'(X \times Y) \simeq \mathcal{D}'(X) \widehat{\bigoplus} \mathcal{D}'(Y)$
$\simeq L(C_{0}^{\infty},\mathcal{D}')$, where $X,Y$ are open sets in $\Rn$. For instance $U : C_{0}^{\infty} \rightarrow \mathcal{D}'(X)$
as an annihilator, has a representation through $\mathcal{D}'(X \times Y)$ (\cite{Treves67}). 

\newtheorem{lem4}[prop2]{Lemma}
\begin{lem4}
 Assume $\mid \eta \mid \leq 1$, is the support for $F (U^{\bot_{E}} \gamma)$.
Thus when the movement changes character, the support is not one-sided in $\eta$. A transmission property, in this context
is a proposition that for instance $U_{1}^{\bot_{E}} =0 \rightarrow U_{2}^{\bot_{E}}=0$ is projective.
\end{lem4}
Consider $F(U^{\bot_{E}} \gamma)=\widehat{V}$ and the mapping $R_{U} : U \gamma \rightarrow \eta$, the reflection axes,
thus the range of $R_{U}(U \gamma)=\{ \eta \quad U^{\bot_{E}} \gamma=0 \}$. When $\widehat{V}$ has compact sub-level sets, 
we consider $\mbox{ sng }\widehat{V}$
as the set $\{ \widehat{V}=\frac{\delta}{\delta \zeta_{j}}\widehat{V}=0 \}$, that is reflection points are singular 
points for $\widehat{V}$, when the axes is not dependent of $\gamma$. When the movement changes character,
the mapping $R_{U}$ has jump discontinuities.  When the two reflection axes are given by $R_{U_{1}}$ and $R_{U_{2}}$, we assume a movement $U_{0}$
that combines the two axes, $U_{0}$ is assumed independent of $U_{1},U_{2}$ and $U_{0} \gamma=\gamma$
through this continuous movement. In particular when $\eta(\zeta) \sim_{0}$ polynomial and $S_{c_{1}}=\{ \eta=c_{1} \}$
a first surface to $\eta$ and axes for invariance. Then (\cite{Nishino68}) we have a continuous path $S_{c_{1}} \rightarrow S_{c_{2}}$.
As $\widehat{V}=0$ when $\eta=c_{j}$, the path is $\subset $ this zero-set and $\frac{\delta \widehat{V}}{\delta \zeta_{j}} \equiv 0$,
that is singular points for $\widehat{V}$.

According to (\cite{Lie91}, chapter 6 Satz 10 and the following example), when $X dy - Y d x=0$ is the differential equation
representing a one-parameter group $Uf=\xi \frac{\delta f}{\delta x} + \eta \frac{\delta f}{\delta y}$, we can using change of variables $(x,y)$ to $(\xi,\eta)$ put the
equation on separated form, why it can be integrated using quadrature. We use particularly $\xi=x$,
$\eta=y/x$ and $y'=\varphi(\frac{y}{x})$, why the equation can be written $\frac{\delta f}{\delta x} + \varphi(\frac{y}{x}) \frac{\delta f}{\delta y}=0$.

Assume $\eta=e^{\phi}$ and $U^{\bot_{E}} \eta \sim e^{W^{\bot_{E}} \phi}$,
then $U^{\bot_{E}} \eta = g e^{W^{\bot_{E}} \phi}$, where $g \neq 0$ $\log g + W^{\bot_{E}} \phi \in L^{1}$ (\cite{Riesz56}, Radon-Nikodym theorem). 
When $W^{\bot_{E}} \phi \in L^{1}$, it has isolated singularities. Note that over $C$, we have that $U \rightarrow {}^{t} U$
preserves character of movement. $e^{W \phi}-e^{\phi}=0$ iff $e^{W \phi} = e^{\phi}$ iff $e^{W \phi - \phi} = 1$
iff $W \phi = \phi$. Geometrically we can write in $L^{1}$, $U^{\bot_{E}} =0$ iff $W^{\bot_{E}} =0$. The singularities in $\zeta$
are isolated for $W^{\bot_{E}} = 0$ and algebraic for $U^{\bot_{E}} =0$. When the movement is composite, $W^{\bot_{E}} \phi$ changes character iff
$U^{\bot_{E}} \eta$ changes character. Note that $\mid \eta \mid < 1$ can be represented $\{ \phi < 0 \}$
For a pseudo-convex function $\phi$ a neighbourhood of $C$ can be represented this way. 
When the movement is composite $U = U_{1}U_{2}$, we consider $\{ \phi_{1} > 0 \} \cup \{ \phi_{2} > 0 \}$.
Since $\mbox{ max} \{ \phi_{1},\phi_{2} \}$ is pseudo-convex, the neighbourhood of $C$ can still be given
by a pseudo-convex function, as above (\cite{Oka60}).

\subsection{Density for kernel to localizer}
When the operator corresponding to $E$ is considered in $Exp$ (\cite{Martineau}), we can use $L^{1}$
norm for the phase to $E$. We will in this section assume the operator corresponding to $E$ is determined
by the behaviour of $E$ relative some relevant norm, such that the operator is Fredholm. For instance when $\gamma$
polynomial, we have $\gamma(D) G \rightarrow E(\gamma)$, where $G$ is Fredholm. As long as we have consecutive
simple movements, this notation is applicable. When we consider spiral movements, the operator is not assumed 
Fredholm.

\newtheorem{prop3}[prop2]{Proposition}
\begin{prop3}
 The Schwartz kernel to the localizer, corresponding to a partially hypoelliptic operator, over 
 the kernel (zero-space), corresponds to a hypoelliptic operator.
\end{prop3}
Consider $E$ as corresponding to a Fredholm operator, with trivial kernel for $E^{2}$. In this case,
$<E(\phi),{}^{t} E(\psi)>=0$ , for all $\phi \in (I)$, implies $\psi=0$.
Note that when $(I- E) \rightarrow (I + E)$ bijective, we have $E^{2}$ corresponds to a very regular action, 
if this holds for $E$.
Note that in this case, $E^{\bot_{E}} \varphi=0$ implies $E(\varphi)=\varphi$, that is $\varphi=0$ (modulo $\widehat{C}^{\infty}$). 
When $E$ is to a hypoelliptic operator,
the mapping $E^{\bot_{E}} \rightarrow (-E)^{\bot_{E}}$ is bijective. When $E^{N}$ is to a hypoelliptic operator,
 we require that the mapping $\overline{E}(\varphi) - \overline{\varphi} \rightarrow E^{\bot_{E}} (\varphi)$ is bijective.

When $<H,E^{\bot_{E}}>=0$, we have $\big[ {}^{t} H , E^{\bot_{E}} \big] \bot_{E} I$ and $\big[ {}^{t} H , E^{\bot_{E}} \big]=0$.
When $E^{\bot_{E}}$ is closed, we have $R(E^{\bot_{E}})={}^{\circ} N(E')$ (annihilator). Assume $E$ with Fredholm
representation and $R(E^{\bot_{E}}) \bot_{E} R(E)$, it remains to prove that the annihilator $H$ is hypoelliptic. Note that
$E^{N}$ is hypoelliptic on $R(E^{\bot_{E}}) \bigoplus R(E)$ and $E^{N} \bot_{E} E^{N-1}$ on $R(E^{\bot_{E}})$. This can be used
to construct $H$.

When $\big[ H,E \big]=0$ and $H$ hypoelliptic, we can conclude $E$ corresponds to regularizing action. 
In the same manner
$\big[ H,E^{\bot_{E}} \big]=0$ implies $E$ corresponds to very regular action. Assume $N(E)=\cup_{1}^{N} X_{j}$, where 
$X_{j}=N(E^{j}) \backslash N(E^{j+1})$ a ``stratification''. Note that when $E^{2}-I$ to regularizing action
and $E-I \rightarrow E + I$ bijective, we can conclude $(E-I)^{2}$ corresponds to regularizing action on $X_{1}$. Note that
when $\big[ E^,I \big] \sim \big[ I,E \big]$, $E^{\bot_{E} 2} \sim E^{2 \bot_{E}}$. We argue recursively, $X_{N}=\{ 0 \}$. Thus, there
is a reduced $H \bot_{E} (E^{N})^{\bot_{E}}$ such that $E^{N}$ corresponds to very regular action on $X_{N-1}$.
When $E^{N-1} \bot_{E} E^{N}$ on $X_{N-1}$,
we can represent $E^{N-1}$ as to very regular action on $X_{N-1}$. 
In the same manner $E^{N-2}$ to very regular action on $X_{N-2}$,
finally $E$ to very regular action on $X_{1}$. 

Assume $E$ an integral operator with kernel $E$, then we write $I_{E}(f)=\big[ E,f \big]$ and $I_{{}^{t} E}(f)=\big[f ,E \big]$.
Consider $<P \big[ E,f \big],\psi>=<P \big[ f,E \big],\psi>$ and assume $P$ such that $P E= E {}^{t} P$ corresponds to 
very regular action
outside the kernel to $E$. In this case modulo regularizing action, $\big[ I,f \big]=\big[ f,I \big]$.
Assume $f =P f'$, where $f'$ analytic and $P$ polynomial, then a sufficient condition to conclude $f$ algebraic
in the sense that $\big[ I,f \big]=\big[ f,I \big]$, is that $\big[ E,f \big]=\big[ f,E \big]$, for 
$E$ such that $P E$ corresponds to very regular action and $ PE=EP$, that is $E$ corresponds to two-sided parametrix. Note that in this case 
$P$ can be used as a base for the ideal, $P = \gamma$. 

When $M$ is localizer to $\eta + i v$, it can be represented as to very regular action.
When $\mid M(\eta,v) \mid \leq C \mid (\eta,v) \mid \leq \lambda$ and denote $W$ for the sub-level sets 
to $M$ and $V$ for the sub-level sets to $\eta + i v$, then we have that the sets $W \cap V$ are compact.

 Assume $K_{j}$ the kernel that represents the movement in simple movements. Assume we have existence of
 $\lim_{j \rightarrow \infty} K_{j}$ in $\mathcal{D}'$. For instance when $K_{j}$ a.c. and the limit
 is taken as $\Sigma (I_{j} - I) \rightarrow 0$, for parameter intervals. Thus the spiral movement, is
 approximated by factorized movements. When we have density for range modulo regularizing action,
 we can use the generalized MP to conclude existence of limit. 

\subsection{The generalized moment problem}
 Schr\"odinger operators represent
a global model, that is not normal. The movement corresponding to adding a potential, is assumed such that
$U_{1}U_{2} \gamma \leq U \gamma \leq U_{2}U_{1} \gamma$. In this case we do not have an algebraic base, but an integrable
base. 

Assume the movement has a factorization or can be approximated by a factorized chain. When $F$ is monotonous
locally, we have $F(W_{1} \gamma) \leq F(U \gamma) \leq F(W_{2} \gamma)$, where $W_{j}$ are the approximating chains.
We can now use the generalized moment problem, that is we consider $\limsup$ and $\liminf$ over factorized chains.
The proposition is that every element in $C^{0}$
can be approximated arbitrarily close by $E_{0}$
iff we have existence of $\alpha(x)$ of bounded variation, such that $\int g(x) d \alpha(x)=0$ 
for all $g \in E_{0}$ implies $\int g(x) d \alpha(x)=0$
$\forall g \in C^{0}$ (\cite{Riesz56}).
Note that when $\widehat{V}$ can be defined relative the factorized chains, the potential to $U \gamma$
can be constructed as a limit.

The measure associated to the movement on the cylindroid can be assumed analytic, why we can construct $d w(z,\overline{z}) \sim \mid d w(z) \mid$
assuming $z \rightarrow \overline{z}$ projective. Consider $\gamma$ restricted to a division element and 
$F$ continuous on $\gamma$ and complex valued,
and $dF$ bounded variation on the division element. For the restriction
to division elements, we can write $<F(U^{\bot} \gamma),\phi>=\int \phi d \Xi(U^{\bot})$, where $\phi = \phi(x) \in (I)$ and
$(I)$ is defined by the movement. Since the spiral movement is formed in the kernel to $F$, using density
for the kernel to the localizer, we can construct the spiral as the limit of factorized movements, as
the division gets finer.

\section{The inverse lifting principle}
We consider $F(\gamma) \rightarrow \gamma \rightarrow \tilde{\gamma} \rightarrow \zeta$, where $\gamma \rightarrow \tilde{\gamma}$
continuation of $\gamma$. The problem is when $\tilde{\gamma}=U \gamma$ can be selected as algebraic. We look for sufficient
conditions on $\widehat{V}$ to determine $U$ uniquely. 
When $F(U \gamma - \gamma)=\widehat{V}=0$ implies $U \gamma =\gamma$, we can conclude axes for invariance and movement.
A sufficient condition for this is that $F$ has a trivial kernel.
Note that $U$ is a functional, why
monodromy is not applicable. 
\subsection{Involution representation}
Consider the model $F(\gamma) \rightarrow \gamma \rightarrow \tilde{\gamma}$, where $\log \tilde{\gamma} \in L^{1}$.
The condition for involution in $(I_{1})$ has an exact correspondent in $(I_{2})$, but the condition for 
orthogonals are related to different metrics (scalar products). 
Assume the zero's to $\frac{d y}{d x}-1$ indicate change of character for the movement. Assuming
$\frac{d y}{d x}$ has isolated singularities, the same hold for $\frac{d y}{d x}-1$ that is
$\{ \frac{d^{2} y}{d x^{2}}=\frac{d y}{d x}=0 \}$ are isolated. 
When $\frac{d y}{d x} \rightarrow 0$ in $\infty$, this means that the movement does not change character
in $\infty$, that is in a neighbourhood of $\infty$, $\frac{d y}{d x}$ is regular and $\neq 1$.

Consider the following examples. When $\widehat{V}(\zeta_{t})=const$ on a connected domain, we have
$\frac{\delta}{\delta t} \widehat{V} \equiv 0$ on the domain, where $t$ is movement parameter. Sufficient
for this is a single movement, but we can not determine the character of movement without conditions on the domain.

When $\widehat{V}=\delta_{0}$, consider $\phi, \frac{d}{d x} \phi \in \mathcal{D}(\Omega)$ and $\phi(0)=1$, that is
$<F(U^{\bot} \gamma) , \phi>=1$. Assume for instance $F=\frac{d}{d x}G$ and $\mbox{ ker }G=\{ 0 \}$,
then $<G(U^{\bot}), \frac{d}{d x} \phi>=<H,\frac{d}{d x} \phi>=1$. Then we can conclude for $U^{\bot} \gamma=0$
implies $\mid x \mid> \mid y \mid$, that is $U=U_{1}$, a single movement. The example can be extended to higher order derivatives.
Further, using monotropy, we can give an argument similar to the previous example, and given conditions
on the domain, get a conclusion for $U_{3} U^{\bot}$, where $U_{3}$ is parabolic.
\newtheorem{def4}[prop2]{Definition}
\begin{def4}
 Involution is defined using the Poisson bracket $\{ F_{1},F_{2} \}=$
 $\Sigma_{i=1}^{p}$ $(\frac{\delta F_{1}}{\delta x_{i}} \frac{\delta F_{2}}{\delta y_{i}} -$ 
$\frac{\delta F_{1}}{\delta y_{i}} \frac{\delta F_{2} }{\delta x_{i}})$. Assume $F_{k}$ for $k=1,2$ satisfies
$\frac{\delta F_{k}}{\delta x_{i}}=-Y_{k}$, $\frac{\delta F_{k}}{\delta y_{i}}=X_{k}$, then
when $p=1$, the condition above is Lie's condition $Y_{2}/X_{2}=Y_{1}/X_{1}$
\end{def4}

The linearity condition $\frac{d y}{d x}=const$, given a contractible domain implies
$\frac{d F}{d x} / \frac{d F}{d y}=const$ on a set of possibly positive measure. The same argument
for $\frac{d F}{d \zeta_{1}} / \frac{d F}{d \zeta_{2}}$ is impossible if $d \zeta_{1} \bot_{E} d \zeta_{2}$.
For instance, when $\frac{d y}{d x}$ constant and $\frac{d F}{d \zeta_{1}}/\frac{d F}{d \zeta_{2}} \rightarrow 0$,
when $\zeta_{1}=P(\zeta_{2},\ldots,\zeta_{n})$ with $\frac{d P}{d \zeta_{2}}$ reduced.
The conclusion is that we may have a linear dependence in $\gamma$ and at the same time 
a non-linear dependence in $\zeta$.

\subsection{A one-sided lift operator}

Assume $T$ is a distribution such that $T F(\gamma)=F(U^{{\bot}_{E}} \gamma)$,
for all $\gamma \in \mathcal{D}$. Over generators, we can assume $T \equiv 0$ over for instance $\mid \eta \mid \geq 1$.
Over the generators (assume $\gamma \in \mathbf{R}^{2}$), $T=\big[ T , H \big]= \big[ H, T \big]$, where $H$ is Heaviside, that is $d T$ is algebraic over generators.

Any movement is assumed to satisfy a local Pfaff condition, when the condition is global
in phase and if the symbol is regular analytic, we have that the movement is simple in the euclidean space. The blow-up mapping $\Phi(\frac{d y}{d x})=\eta$ is thus dependent on the movement,
when $\Phi'=\Phi < 1 (=1, > 1)$ the movement is translation (scaling, rotation). Note that the involutive condition in $V$, means that $J$ preserves
symmetry in $(x,y)$. Assume $\Sigma=\{ (x,y) \quad U_{1} \gamma = U_{2} \gamma \}$, thus $L U_{1} = L U_{2}$
and $J \Sigma =\{ E V_{1}  = E V_{2}  \}$, thus a symmetry set in Lorentz metrics has a correspondent 
in euclidean metrics. Assume $\tilde{F} \sim J F$, then ${}^{t} F \sim F$ over $U$ has a correspondent relation
${}^{t} \tilde{F} \sim \tilde{F}$ over $V$. 

If we assume the movement monotonous, for instance such that the (euclidean) distance to axes $\eta=c$
is increasing, as the parameter grows, the support for $\widehat{V}$ is one-sided in the plane.
We assume when we have an algebraic base $\gamma$, that we have a transmission property. Simple movements
preserve this property, that is when $C$ is discrete, we still have a transmission property. When the kernel
to $M$ is symmetric and such that $(\eta,v) \rightarrow (v,\eta)$ is projective,  in case $\eta^{2}$ reduced,
$\eta=v$ is discrete. Note that over $C$ we consider $\tilde{U} \eta \sim \mid U \eta \mid$, why we must
require $\eta \overline{\eta}$ reduced. Note that in the spiral case, we do not assume $\eta \rightarrow v$ 
projective in limit.

\subsection{Projectivity for base space}

Assume $\psi(\eta)=\frac{1}{\eta}$ and $U$ locally algebraic, with $M(U \eta,v)=M(\eta, U v)$.
Then when ${}^{t} \psi = \psi$, we have ${}^{t}M = M$. When $\eta^{2}$ reduced, we have $\eta = v$ implies $\eta=\eta_{0}$,
a point. Thus the domain for $F$ can be written $\eta \bigoplus v$ and $\psi$ is projective. 

 Given $U F \sim F U$, when $\zeta_{t}$ 1-1 and closed, then there
is a connected $\{ U \gamma \}$ such that $T F(\gamma)=F(U \gamma)(\zeta)=F(\gamma)(\zeta_{t})$.
Otherwise we have that $\zeta_{t}$ is dependent on $U \rightarrow T$ or
$F U - U F$. In particular, we have existence of $U \gamma$ such that $F(U \gamma)(\zeta)=f(\zeta) + \lambda \delta_{0}$.
When $\zeta_{t}$ is not 1-1, we have that $U \gamma \in X \bigoplus X_{0}$, where $U \gamma$
is in-determined on $X_{0}$, that is $f(\zeta) + \lambda \delta_{0} = F(U \gamma + U \gamma_{0})$.
Alternatively, consider $T F(\gamma)=F(U \gamma) + \widehat{V}$
where $\widehat{V}$ has support corresponding to the set for indetermination.

Assume $(IP)=\{ \gamma \in H_{m} \quad \mid \eta \mid < \lambda \quad \eta \mbox{ polynomial } \}$. We consider $(I) \sim (IP) \bigoplus (I_{2})$ and
$(J)=\{ \eta \quad V^{\bot_{E}}  \eta = 0 \}$. Thus when $(J) \subset (IP)$ we have a single movement. When $\eta \rightarrow 1/\eta$
is projective, we can assume $(I_{2})=\{ 1/v \quad v \in (IP) \}$. 

\subsection{Composition of movements}
Assume $A : \mathcal{D} \rightarrow \mathcal{D}$, then we can obviously define $\big[ F,A \big] \in \mathcal{D}'$,
when $F \in \mathcal{D}'$. When $A : \mathcal{D}_{L^{1}} \rightarrow \mathcal{D}_{L^{1}}$ and we have a very regular boundary (\cite{Dahn13}), 
we can relate the functional to a polynomial. 
A very regular boundary, can be represented as $\{ F_{T}=const \quad \frac{d F}{d T}=const \}$ where
$F_{T}$ is holomorphic otherwise.

When $I$ denotes a mapping invariant over an axes, we must have for composite movements $\big[ I,A \big]=\big[ A,I \big]$,
that is the movement as composite with reflections is assumed algebraic. When $A,B,C$
are reflections, $A (BC)=(AB)C$. For simple movements, we have that $A \gamma$ is polynomial, if $\gamma$ is polynomial.
Over a cylindrical domain, we assume
$2 \log I \sim \log I^{\bot_{E}} + \log (-I)^{\bot_{E}}$. Further when $U_{1} \bot_{E} U_{2}$ with respect to a two-dimensional (projective)
space, with identity $J$, we assume $2 \log I \sim \log \overline{J}^{\bot_{E}} + \log (-J)^{\bot_{E}}$.

When $U_{1}=U_{2}$ and $\frac{d U_{1}}{d U_{2}}=1$ we have $\gamma=\gamma_{0}$ that is a point.
Assume $U_{1},U_{2}$ have compact sub-level sets, then the sequential composition $U_{1} U_{2}$ has compact
sub-level sets. Conversely, when the $U \sim U_{1} U_{2}$ has compact sub-level sets and $ U \gamma (\zeta) \rightarrow \mid \zeta \mid$
is continuous, these compact sets have a division of $\mid \zeta \mid$. When $J$ is proper,
the division has a corresponding division of $V \sim V_{1} V_{2}$, as long as $V \gamma (\zeta) \rightarrow \mid \zeta \mid$
continuous.

Assume $M \sim {}^{t} M$ and $\mbox{ ker }M $ non-trivial, where $(\eta,v)$ has compact sub-level sets.
When $\eta^{2}$ reduced, the mapping $\eta \rightarrow v$ is projective. When
$(\eta,v)$ describes a full line $\subset \tilde{C}$ (limit of $C$), we note that $H_{m} \cap \tilde{C} \rightarrow \zeta_{t}$ is continuous.
Outside $\tilde{C} \cap H_{m}$ $(\eta,v)$ may still have compact sub-level sets, but we do not have a continuous inverse,
$U \eta (\zeta) \rightarrow \zeta$.

\section{Dependence on singularities}
\subsection{Critical points}

When the movements are considered as holomorphic mappings, the set $C$ is critical for some
changes of character of movement and not for others.
Since the movements are considered as functionals, the set $C$ is not considered as critical
for the movement. However, when $\widehat{V}$ is considered as regular in $\zeta$, the pre-image to $C$ is considered
as singular for $\widehat{V}$.

When $\eta$ is analytic, any pole $a$ is isolated and it is not a pole for $v$.  
Given that $F$ is analytic, such that we have existence of $P,Q$ polynomial with  $F=P/Q$, when $Q=const$, $F$ is algebraic. When $Q$ reduced, $\{ Q=const \}$ is algebraic.
When further $P$ is invertible (outside zero's for $P$ polynomial), we have $F/P \sim 1/Q \rightarrow 0$ in $\infty$

Assume $\log y \in L^{1}$. On the set where $\frac{d y}{d x}=1 \notin L^{1}$, but analytic, we do not have isolated singularities
but algebraic singularities. If the set is considered as boundary to $\frac{d y}{d x} \geq 1$, in particular $y$ can be considered
as monotonous. 

The condition $\frac{d y'}{d x'}=\frac{d y}{d x}$ analytic means that $reg \rightarrow reg$.
The corresponding parabolic condition, given it defines an analytic set means $\eta=1 \rightarrow \eta'=1$
continuous. 
Given that $\gamma$ hypoelliptic, that is $x \prec \prec y$, we have that $\{ \zeta \quad \frac{d y}{d x}=1 \} \subset \subset \Omega$
that is ``jumps'' do not contribute with singularities micro-locally.

Assume $f(\zeta_{t})= A_{t} g(\zeta) = (1 + \widehat{V} / g)$. When $g$ is considered over an a.c. domain
$\frac{\delta }{\delta \zeta_{j}} g=0$ implies $g = const$. When $\frac{\delta f}{\delta t} = \frac{\delta}{\delta t} A_{T} g(\zeta) \sim \frac{1}{g} \frac{d \widehat{V}}{d t}$,
we assume that $\gamma$ is dependent on $f,F$ and $U$ is dependent on $\zeta_{t}$ and $\widehat{V}$.

\subsection{Orientation}

Assume the boundary $C$ is defined by $\varphi(\zeta)$, when a neighbourhood of $C$ is pseudo-convex
and $\varphi$ pseudo-convex, the neighbourhood can be written $\{ \varphi < 0 \}$. A necessary condition
for strict pseudo-convexity, is that the normal has a local algebraic representation (\cite{Oka60}). A necessary condition for this
is that the boundary is oriented (\cite{Dahn13}). Assume $\varphi^{*}(\zeta)=\varphi(R \zeta)$, where $R \zeta$ is planar. 
When $\varphi^{*}=-i \varphi$,
analyticity is preserved.

For a composite movement, the regularity is dependent on orientation with respect to SGN function.
When the movement can be divided into single movements, we consider it as monotonous on segments.
In the case where the movement changes character on a segment, we assume it does not simultaneously change
orientation. 
The $SGN(p_{1},p_{2})$ function orders the zeros to two real polynomials in one real variable,
and the intervals between them (including $\pm \infty$) (\cite{Garding97}). In particular when $\log U^{\bot_{E}} y \in L^{1}$
we assume $U^{\bot_{E}} y= y(U^{\bot_{E}})=y(u_{1},u_{2},u_{3})$, geometrically equivalent with a polynomial 
in parameters. When $\eta'=const$ an axes, we have $\frac{d}{d x'} (y' - x')=0$. When the movement changes character,
using involution and geometric ideals,
we can write $U_{j}^{\bot_{E}} y = A_{j} \frac{d y}{d x} + r_{j}$, for regular $A_{j},r_{j}$. Using the collar point condition,
when $\frac{d}{d t} U_{1} r=\frac{d r}{d u_{1}}$, we have $\frac{d r_{1}}{d r_{2}}=\alpha$.
Given that $U(t) y$ has a semi-algebraic sub-graph, the movement can be divided into finitely many locally
analytic, algebraic movements (\cite{Garding97})

Assume $\gamma$ polynomial and consider $(e^{v},e^{\phi + v})$, where $\eta=e^{\phi}$ and the movement
is defined by $U \eta \sim_{0} e^{W \phi}$. Thus, we are assuming a global definition $W_{n} \phi \in L^{1}$,
where $W_{n}$ is dependent on division of parameters. The division is assumed to get finer as $n \rightarrow \infty$.
When the movement changes character, we are assuming collar points for finite $n$, that is $d U_{1} / d U_{2} = \alpha_{n}$
where $\alpha_{n}$ regular and dependent on division. When $n \rightarrow \infty$, there are cluster sets for the 
sub-level sets to $\alpha_{n}$, but the division in $\mid \zeta \mid \rightarrow$ a point, as $n \rightarrow \infty$. Let $\phi_{n}=W_{n} \phi$ and assume $\parallel \phi_{n} - \phi_{m} \parallel \rightarrow 0$
as $n,m \rightarrow \infty$ and $\phi_{n} \rightarrow \phi$ in mean convergence. Thus, we have existence
of $\lim_{n \rightarrow \infty} U_{n} \eta$ in $Exp$ and $\mathcal{D}_{L^{1}}'$, as $\phi_{n} \rightarrow \phi$
in mean relative $L^{1}$.

Assume $\phi_{n}=\phi(W_{n})$, then as $W_{n} \rightarrow W$, we have $\phi(W_{n}) \rightarrow \phi(W)$
continuous. However, as $n \rightarrow \infty$, the linear independence for parameters is not preserved.
When $\phi_{n}(W) \equiv \phi_{n}(W')$, for $W \neq W'$, then $\phi_{n}$ is a point for finite $n$.
As $n \rightarrow \infty$, $\phi_{n}$ can be connected. When $\parallel \phi_{n} \parallel \rightarrow \parallel \phi \parallel$
implies $\phi_{n} \rightarrow \phi$ (strong convergence), the limit movement $W$ can be given a global
definition. When $\eta \rightarrow v$ is projective in $W_{n}$ as $n \rightarrow \infty$, we have a global definition
of $W$ as $n \rightarrow \infty$. In this case there is no room for a spiral.

\subsection{Completion}

Consider $(F (\gamma), M(\eta,1/\eta))$ and $F (U \gamma)=F(\gamma) + \widehat{V}$.
Further, $({}^{t} F(U \gamma), M(U 1/\eta,\eta))$ where ${}^{t} F(U \gamma) = {}^{t} F(\gamma) + {}^{t} \widehat{V}$.
When $T F(\gamma)(\zeta)=F(U \gamma)(\zeta)$ analytic, for instance over $C$, we define the inverse $\bot T F(\zeta) = \zeta_{t}$,
where $t$ is parameter for a single movement.

\newtheorem{lem3}[prop2]{Lemma}
\begin{lem3}
Given that $F(U \gamma)$ analytic in $\zeta$ and $\bot T F(\gamma)$ is locally 1-1 and closed, we have that 
$T$ is surjective, that is $T f(\zeta) = f(\zeta_{t})=F(U \gamma)(\zeta)$
for some $\gamma$. Otherwise (continuum), the system $\{ U \gamma \}$ must have a non trivial
complement where $\widehat{V}_{t}$ completes the model.
Thus $f(\zeta_{t}) - f(\zeta)=F(U \gamma) - F(\gamma)=\widehat{V}(\zeta_{t})$ $(= F(U \gamma - \gamma)$ when $F$ linear).
When $\widehat{V}$ is given so that $\zeta_{t}$ is closed and locally 1-1, we have existence of $\gamma_{t}$
such that $F(\gamma_{t})=F(\gamma) + \widehat{V}_{t}$.
\end{lem3}

(\cite{Martineau} Ch I, Theorem 4)
In this case $(U - 1)\gamma \sim F^{-1} \widehat{V}$
dependent on involution. We assume $\widehat{V} \equiv 0$ on $\Delta$ (lineality). 
When $F^{-1} \widehat{V}$ is regular, the movement
$U \gamma - \gamma$ can be determined. Assume $\gamma_{T}^{N}$ reduced with localizer $F_{N}$ 
and  $\widehat{V}_{N}$ the corresponding potential, such that $\widehat{V}_{N}$ is locally 1-1.
Assume $N(\widehat{V})=$ lineality. If $F^{2}$ reduced, $\widehat{V}_{2}=0$ implies $\zeta=0$ and 
$\widehat{V}_{2}=2 F \widehat{V} + \widehat{V}^{2}$. Concerning the involution, when $d F=0$, on
a contractible domain $\frac{d y}{d x} \sim \frac{d F}{d x}/\frac{d F}{d y}$.
$\frac{d F}{d \zeta_{j}}=(-Y + X \frac{d y}{d x}) \frac{d x}{d \zeta_{j}}$. Thus on a contractible domain,
when $F$ is Hamilton function, $\frac{d F}{d \zeta_{j}} \equiv 0$. The involution condition is such that
$\frac{d F_{T}}{d \zeta_{j}}=0$ implies $\frac{d F}{d \zeta_{j}} =0$.

Consider $f=F(\gamma)=F(U \gamma + U^{\bot_{E}} \gamma)=F(U \gamma) + \widehat{V}$. 
Consider $U^{\bot_{E}}$ as an extension mapping.
The movement $U^{\bot_{E}}$
is always defined in Banach-spaces, for instance Exp. 
The mapping is not involutive in the entire space, but if we consider a cylindrical domain of movements $\tilde{U}$,
such that $(I_{3}) \sim (I_{1}) \bigoplus (I_{2})$ and when $(I_{1})^{\bot_{E}} \bigoplus (I_{2})^{\bot_{E}} \sim (I_{3})$,
then it is considered as involutive over this boundary. 

Consider $F(\gamma) \rightarrow \gamma \rightarrow \tilde{\gamma} \rightarrow \zeta$ and 
$F(\tilde{\gamma})=F(\gamma) + \widehat{V}$, that is $\widehat{V}$ is defined relative the continuation of
the domain for $\gamma$. 
We can construct $\tilde{\gamma}$ using standard movements $U_{j}$, that is $U_{1} U_{2} \gamma \leq \tilde{\gamma} \leq U_{2} U_{1} \gamma$.
When $\widehat{V}$ defines a standard movement, $\widehat{V}=0$ implies $\mid \eta \mid < 1$ ($> 1, =1$)
not dependent on point. For the proposition that $\tilde{\gamma}$ polynomial,  when we extend the domain with standard movements, we assume 
algebraicity is preserved.

\subsection{The blow-up mapping}
The condition $x d x - y d y=0$, that is $\frac{y}{x}=\frac{d y}{d x}$ and $\Phi(\frac{d y}{d x})=\frac{d y}{d x}$
means a linear dependence, a ``flat'' blow-up. Let $\psi(x)=\frac{1}{x}$ involution. When
$\psi$ is bounded on $y/x \leq 1$, the movement must change character. Let $\psi_{1}(x-1)=\frac{1}{x}-1$,
that is $\psi_{1}(\frac{y}{x}-1)=\frac{x}{y}-1$. When $\psi_{1}$ is unbounded, we have that the movement
preserves character.
The condition on ``collar points'' means that  $\frac{d U_{1}}{d U_{2}}=1$ can only occur in points, when the movement 
is on the hyperboloid. In non-hyperbolic space,
the axes for reflection, is determined by a condition on $\eta$, that is  the sign of $\Phi-1$ determines the axes.

 Under the movement, we have either that the axes is fixed $\frac{d y}{d x}=const$
or moving $\frac{d y}{d x}=\rho$ regular. The involution condition $\frac{d y}{d x}=\rho$, for moving axes, 
is assumed such that $\rho \rightarrow 1$ regularly, that is using a scaling of hyperboloid, we can assume change
of character occurs for $\rho=1$.
Assume $W$ is such that $\rho=1$. 
The involution does not completely determine the movement, consider $\Phi(\frac{d y}{d x})=\eta$, that is
$x {}^{t} \frac{d}{d x}\Phi(y)=y$. Thus if $x < y$ then ${}^{t} \frac{d}{d x} \Phi > 1$, when $x=y$ 
then ${}^{t} \frac{d}{d x} \Phi=1$. If $y=\frac{d z}{d x}$ and
${}^{t} \frac{d^{2}}{d x^{2}} \Phi > 0$, then for strong derivatives $\Phi(z)$ is convex.

Lie (\cite{Lie91})considers $\frac{d y}{d x}=\varphi(\frac{y}{x})$, that is $\Phi \varphi=id$. Assume $\xi=\frac{y}{x}$ and $\eta = x$.
 Thus when $\frac{\delta f}{\delta x}=-Y$,$\frac{\delta f}{\delta y}=X$, then 
$\frac{d y}{d x}=\frac{Y}{X}=\varphi(\frac{y}{x})$. 
According to Bendixson (\cite{Bendixsson01}), when $(X,Y)$ analytic determines a dynamical system with
trajectories $\Gamma$, then $\Gamma \rightarrow 0$ implies $\Gamma$ a spiral or $x Y - y X =0$
(determined tangent), that is $\Phi \equiv I$.

\newtheorem{def2}[prop2]{Definition}
\begin{def2}
 Assume $C$ the set where the movement $U$ changes character in a cylindroid movement,
 that is $\log U_{1} \sim \log \tilde{U}_{2}^{\bot_{E}}$ with a collar point according to $d U_{1} / d U_{2}=const$
and $U_{1}=U_{2}$ implies $\gamma=\gamma_{0}$, a point.
\end{def2}

The set $C$ are not singularities for the movement relative $x$, but singular for the correspondent 
movement in $\zeta$. A point can be considered as isolated, given $U^{\bot_{E}}_{1} x_{n} \rightarrow 0$, $U_{2}^{\bot_{E}} x_{n} \rightarrow 0$
and $d (U_{1} x_{n} - U_{2} x_{n}) \rightarrow 0$, as $x_{n} \rightarrow x_{0}$, that is $x_{0}$ is a point
(collar point). We can put compatibility conditions such that
$U_{1} \eta - U_{2} \eta \rightarrow 0$ iff $U_{2}^{-1} U_{1} \eta - \eta \rightarrow 0$
and $d (U_{2}^{-1} U_{1} \eta - \eta) \rightarrow 0$, implies $\eta = \eta_{0}$, a point.
The proposition on isolated singularities, can be put as $\zeta$ regular for $U^{\bot_{E}} x_{n}$, for all $U$
implies $\zeta$ regular for $x_{n}$. When singularities are isolated, $\zeta$ is singular for $x_{n}$ if 
$\zeta$ singular for $U x_{n}$ implies $U = I$.

\newtheorem{prop1}[prop2]{Proposition}
\begin{prop1}
 Given $J$ continuous maps collar points in $U$ on collar points in $V$,
we can represent $J$ as a proper mapping.
\end{prop1}

Assume the movement generated by a parameter $t$, such that $\frac{d U_{1}}{d t}=\alpha \frac{d U_{2}}{d t}$
where $\alpha$ regular outside a discrete set (cf very regular boundary). Note that $\alpha$ is not $0$
when $t \rightarrow \infty$, since this would imply inclusion between corresponding ideals. Using the Radon-Nikodym theorem,
when $\Omega_{2}$ is a neighbourhood generated by $U_{2}$ and correspondingly for $U_{1}$, we have $\int_{\Omega_{1}} \gamma d U_{1}  \sim 
\int_{\Omega_{2}} \gamma \alpha d U_{2}$, where $\alpha \in L^{1}(d U_{2})$. Note that when $\Omega \cap C \neq \emptyset$,
we consider divisions $\Omega = \Omega_{1} \cup \Omega_{2}$, with boundary $C$. We assume $\infty \notin  C$.
We require from $J$ that $\frac{d}{d t} J U_{j}=J \frac{d}{d t} U_{j}$. In particular this means that Lie's condition
$\frac{Y'}{X'}=\frac{Y}{X}$ is preserved by $J$ into euclidean metrics.

 Assume $x_{j} (=U x) \rightarrow x_{0}$ then $x_{j}' (=V x) \rightarrow x_{0}'$
and $x_{0} \sim x_{0}'$ (conjugated). According to Riesz-Fischer theorem, $\parallel x_{n} - x_{m} \parallel^{2} \rightarrow 0$
as $n,m \rightarrow \infty$, means that $x_{n}$ has a limit in the mean. Frechet topology
is closed for this convergence.  
We consider $B x_{n} \rightarrow x_{0}$ weakly and $\mid A x_{n} \mid \leq B \mid x_{n} \mid$, 
which implies existence of $g(x)$ measurable with respect to $B$, such that $A (x_{n})=B(g x_{n})$ (\cite{Riesz56}, Ch. 3, Section 63)
when $x_{n}$ summable with respect to $B$, ($A=\tilde{U},B=\tilde{V}$). Note that $x_{n}' \rightarrow x_{n}$ is proper and $x_{n} \rightarrow x_{n}'$
is proper, when we have convergence in the mean.  
When we consider $\tilde{U} \sim \mid U \mid$, that is $\tilde{U}(x) \sim \mid U(x) \mid$, we have
that convergence is convergence in the mean. The separation property is relative $L^{1}$. The moment problem
solvable for $\tilde{U}$ does not imply solvable for $U$.

When $\Phi(\frac{d y}{d x})=\eta$, the case with $\Phi(\frac{d y}{d x})=\frac{d y}{d x}$, implies a determined
tangent, that is we do not have a spiral. Degenerate points for dynamical systems are given by
$\frac{d y}{d x} + \frac{x}{y}=0$ and then the mapping $\eta \rightarrow v$ does not give transmission over
dynamical systems.
Note that when ${}^{t} \frac{d}{ dx} \Phi=const$, 
when $\int_{\Omega} d \Phi=0$, we have $\int {}^{t} \frac{d}{d x} \Phi(y) d x=\int_{\Omega} const d x=0$ that is $m \Omega =0$

Better, when $\eta=e^{\phi}$, we consider $\int_{(\gamma)} \phi d U$ for instance. With this representation $U_{j} \eta \sim e^{W_{j} \phi}$,
where $W_{j} \phi \in L^{1}$. Note that $(\gamma_{1})$ is a translation domain and $(\gamma_{2})$ a rotation domain
and a point inner to $\gamma_{1}$ can be outer to $\gamma_{2}$.
Assume  $\eta=e^{\phi}$ with $\phi$ pseudo-convex, then $\{ \phi=0 \} \sim \{ \eta = 1 \}$, when $\{ \phi > 0 \}$ implies $\eta > 1$
where $\phi=\phi(x)$ gives a one sided neighbourhood, that is the movement does not change character.
Symmetry means that the movement changes character.
When $\{ \frac{d}{d x} \phi=\phi =0 \} \sim \{ \frac{d \eta}{d x} = 0 \quad \eta=1 \}$.

\subsection{Space of movements}

Consider the space of movements, generated by three axes in $\mathbf{R}^{3}$. The condition that the movement can be factorized in base movements, 
implies existence of Schwartz kernel.  We will assume that the movement can be factorized in a chain such that
$\big[ \big[ A,B \big], C \big]=\big[ A, \big[ B,C \big] \big]$, which is possible when the movements are algebraic.
Note that $AB=BA$ does not imply algebraic movements.

Given for instance $U_{1}U_{2}$ (convex) and $U_{2}U_{1}$ (concave), we can find a trajectory between the two movements
as a geometric mean of the factorized movements. In this case given that
we have existence of limes for the factorized movements, we can determine the limes for the mean.

Note that
$\dim U_{1}=1$,$\dim U_{1}^{\bot_{E}}=2$ and $\dim U_{2}=2$,$\dim U_{2}^{\bot_{E}}=1$ and $\dim U_{3}=1$,$\dim U_{3}^{\bot_{E}}=1$
thus when $\dim U$ is not constant or $\dim U \rightarrow \dim U^{\bot_{E}}$ constant, we have change of character.
Consider $U_{1} U_{2}=I$. Joint invariant sets implies invariant sets
 for the composite movement. Movements over the hyperboloid have disjoint invariant sets, however the mapping
 $U \rightarrow V$ can contribute to invariant sets. Note that joint invariant sets implies a linear
 relation between movements, given a positive measure for the invariant sets.
 
  \newtheorem{lem6}[prop2]{Lemma}
 \begin{lem6}
  Consider $\eta \in \mathcal{D}_{L^{1}}'(\mathcal{U})$, where $\mathcal{U}$ is parameter space for
the simple independent movements. In this case we can write $\eta \sim \Sigma D_{U}^{\alpha} \eta_{\alpha}$
where $\eta_{\alpha} \in L^{1}$. Assume $W$ a movement defined relative a given set $C$, such that for all 
$\phi \bot f$, we have $W \phi \bot \widehat{V}$, this is sufficient to determine $U$.
 \end{lem6}

 Assume $\phi$ such that $<F(\gamma),\phi>=0$ and consider $<\widehat{V},W \phi>=0$, that is we consider
$\widehat{V} \rightarrow W$ and by taking inverse $W \rightarrow U^{\bot_{E}}$. Using an euclidean scalar product,
we can assume $\tilde{U}_{1} \sim \tilde{U}_{2}^{\bot_{E}}$, where the movement changes character. Assume the movement $W=W(t)$
dependent on a parameter. When $\mbox{ ker }F$ trivial, for fixed $\gamma,\phi$, this condition determines $U$.
That is, given a trivial kernel to $F$, or when the kernel has a dense (hypoelliptic) representation, 
every movement of base, has a correspondent movement of orthogonal
to range and by using collar points for $W$, we can determine $W$ as $\sim U^{\bot_{E} \bot_{E}}$. Assume
$<\widehat{V},W \phi>=<W \phi \times U^{\bot_{E}} \gamma, {}^{t} F>$, thus if we have an bijection $\gamma \rightarrow \phi$,
so that every movement in $\gamma$ has a unique correspondent movement in $\phi$, the movement $U$ can be determined
from the movement $W$, when we assume reflexivity for $U \rightarrow U^{\bot_{E}}$ over $C$. When $U \phi \in \mbox{ ker }{}^{t} M$,
and $M$ is taken in $\eta,v$, we can define a spiral movement

\subsection{Determination of movement}
 Assume $U \rightarrow I$ exists over $\eta$, when $\lim_{U \rightarrow I} M(U \eta,v)=\lim_{U^{\bot_{E}} \rightarrow I} M(\eta, U^{\bot_{E}} v)$
 we have a two-sided limit. That is if $U \eta - \eta=0$ implies $\mid \eta \mid < 1$, then $\mid v \mid > 1$,
 that is ${}^{t} U v - v=0$ corresponds to ${}^{t} U = U^{\bot_{E}}$. Assume $T M (\eta,v)=M(U \eta,v)$ and ${}^{t} T M(\eta,v)=M(\eta,U^{\bot_{E}} v)$,
 we have ${}^{t} T M={}^{t} M U^{\bot_{E}}$ and when ${}^{t} M=M$, we have outside the kernel to $M$, $U=I$.
 In case of change
 of character of movement, the dimension for $U^{\bot_{E}}$ is affected (parabolic case excepted). 
 When the movement changes character,
 $U_{1}^{\bot_{E}} \sim$ planar movements, on $C$ we are referring to rotations and $U_{2}^{\bot_{E}}$
 is translation, that is the orthogonal to a spiral movement can be a spiral movement.

 Schwartz kernel theorem means that $<I_{K}(\phi),\varphi>=<K,\varphi \times \phi>$. Not all mappings
 $L_{x}^{2} \rightarrow L_{y}^{2}$ have a kernel in $L^{2}$, for instance symmetry in $\mathcal{D}'$ does not imply
 symmetry in $L^{2}$ (\cite{Treves67}). We have that a desingularization, is not necessary 
 in order to determine the movement uniquely. Assume $\Omega_{U}=\{ \zeta \quad M (U -I) =0 \}$. Given $U$ analytic, we have that
 $N(U-I) \rightarrow \Omega_{U}$ is continuous. When $U \rightarrow T$ is analytic, $\Omega_{U} \rightarrow \Omega_{\bot T}$
 is continuous.

When $F$ is not symmetric, we must use two potentials $\widehat{V}$ and ${}^{t} \widehat{V}$. If the potentials
can be defined as independent on $\eta,v$, they are symmetric relative $(\eta,v)$. Assume $\widehat{V}=0$ implies $\mid \eta \mid \leq 1$
and ${}^{t} \widehat{V}=0$, then $\mid \eta \mid \geq 1$, and $\widehat{V}= {}^{t} \widehat{V}=0$ implies $\mid \eta \mid=1$.
When $F$ is polynomial, we have a transmission property for $M$. Assume $M(\eta,v) \sim P(v,\eta)$ where $P$ polynomial,
this means that we have preservation of value in both variables. In this case $C$ is (semi- ) algebraic
and the movement does not change character. When $\frac{\delta}{\delta \zeta_{j}} \widehat{V}(\zeta) = \Sigma (\frac{\delta M}{\delta \eta} \frac{\delta \eta}{\delta \zeta_{j}}
+ \frac{\delta M}{\delta v} \frac{\delta v}{\delta \zeta_{j}})$ and a symmetry condition can be given by
$\frac{\delta M}{\delta \eta} / \frac{\delta M}{\delta v} = - \frac{\delta v}{\delta \zeta_{j}} / \frac{\delta \eta}{\delta \zeta_{j}}$
over a contractible domain.  
When $\big[ M,U^{\bot_{E}} \big]$ is nuclear over $(\eta,v)$, we can write $<\widehat{V},v \times \eta>$.
In particular $\widehat{V}$ is nuclear over $C$, which after a scaling of the hyperboloid, can be 
represented $\eta=v$.

Consider $U \eta = \eta(\zeta_{t})$, where $\zeta_{t}$ is defined by the movement and where we assume
$\eta$ differentiable in $\zeta_{t}$. We write $\frac{d \eta}{d U}=\frac{d \eta}{d \zeta_{t}} \frac{d \zeta_{t}}{d t}$.
Further, $\frac{d \widehat{V}}{d U}=\frac{d M}{d \eta} \frac{d \eta}{d U} + \frac{d M}{d v} \frac{d v}{d U}$.
Note that when the movement is considered over $H_{m}$, the mapping $U \eta (\zeta) \rightarrow \zeta_{t}$
is continuous. Assume $\frac{d \zeta_{t}}{d t}$ bounded. Then, where $\eta$ reduced,
$\frac{d}{d \zeta_{t}} \log \eta \rightarrow 0$ in $\infty$. When $\eta$ is polynomial, $\frac{d}{d \zeta_{t}} \log \eta \leq C$
in $\infty$. 

Assume $F(U^{\bot_{E}} \gamma)=\widehat{V}$ and $\tilde{F}$ the localizer corresponding to $\tilde{F}(V^{\bot_{E}} J \gamma)$.
Define $\gamma^{2}=(x^{2},y^{2})$ and $\tilde{F} J \gamma \sim \tilde{F}_{2} J \gamma^{2}$ and so on.
Note that modulo $\widehat{C}^{\infty}$, $\tilde{F}_{2} \sim \tilde{F}^{2}$. Using the property that $\tilde{F}_{N}$
has a trivial kernel modulo $\widehat{C}^{\infty}$, we can define $(V^{\bot_{E}} \gamma)^{N} \sim W_{1}^{\bot_{E}} \gamma^{N}$.
When $\big[ V,I \big] = \big[ I, V \big]$, and when $V \rightarrow {}^{t} V$ preserves character of movement,
we can assume $W_{1}^{\bot_{E}} = 0$ iff $V^{\bot_{E}} =0$, that is character is preserved. With these conditions,
we can come to a conclusion concerning which movement $\widehat{V}$ corresponds to. Note that the condition
on the set $C$, that $\tilde{U}_{1} \sim \tilde{U}_{2}^{\bot_{E}}$ for instance, means that $U \rightarrow {}^{t}U$
preserves character over the set $C$.
When $C$ corresponds to consecutive simple
movements, define $\tilde{C}$ as the limit over symmetric $(\eta,v)$, considered in $\mathcal{D}'$.
We do not have a continuous mapping $\tilde{C} \rightarrow \zeta$ in this case. The properties of $\tilde{C}$
are dependent on projectivity for $\eta \rightarrow v$. 

 Note that in hyperbolic metrics, we have an order relation as follows. If 
 $\overrightarrow{x} + \chi \overrightarrow{l}=\overrightarrow{x}'$
 and $L(x,y)=c$, then $L(x',y')=\chi c$, where $L$ is the Lorentz metrics. Thus, 
 when $U_{T}$ is a simple movement and $U_{T_{2}}=U_{T_{1}} + \chi$,  
 in this sense $U_{T_{1}} \leq U_{T_{2}}$ is well defined.
  Note that when ${}^{t} M \psi =0$
 implies ${}^{t}(M U) \psi=0$, this property is characteristic for the movement. Assume $(\eta,v) \bot_{E} {}^{t} (M U_{1} U_{2})$
 and $(\eta,v) \bot_{E} {}^{t} (M U_{2} U_{1}) \psi$ and $U_{1} U_{2} \leq U \leq U_{2} U_{1}$. Assume
 further $(\eta,v) \in \Omega$, a domain for ${}^{t} (M U) \psi \bot_{E} (\eta,v)$, for all $\psi \in \mbox{ ker } {}^{t} M$.
 Then we can extend $\Omega$ to the limit when the division gets smaller (assuming separation property).

\subsection{Inclusion of ideals}

\newtheorem{prop4}[prop2]{Proposition}
\begin{prop4}
 Consider the ideals in $\mathcal{D}_{L^{1}}$ and $(I_{1}) \subset (I_{0}) \subset (I_{2})^{\bot_{E}}$ 
where $\rho_{j}$ are the weights to the corresponding  ideals. A sufficient condition for inclusion is
$\rho_{2}/\rho_{0} \rightarrow 0$, $\rho_{0} / \rho_{1} \rightarrow 0$ and $\rho_{2}/\rho_{1} \rightarrow 0$.
Using the collar point condition, we have $\frac{d U_{1}}{d U_{2}} =\alpha$ and $\frac{d U_{1}}{d U_{3}}=\beta$,
where $\alpha,\beta > 0$ regular. These can serve as weights for local ideals, of functions integrable with
respect to movement parameters.
\end{prop4}
In particular, we can find a regular function $\delta$ such that $\frac{1}{\alpha}/\delta \rightarrow 0$
and $\delta / \alpha \rightarrow 0$, which motivates existence of movement $U$, such that $U_{1}U_{2} \leq U \leq U_{2}U_{1}$.

For composite simple movements, for instance
$U_{1}$ the sequential of $U_{2}$. When $\alpha > 0$, the relation is monotonous, that is when
$U_{1}=U_{2}$, the set $C$ is discrete. 
When $C$ is not discrete, we consider
a distributional representation of movements, over the argument $\eta + iv$ where $v \sim 1/\eta$.
We are assuming $d U_{1}(\eta + i v)/d U_{2}(\eta + i v)=\alpha_{n}$ regular (non-constant). 

Assume $d \eta \rightarrow - \eta^{2} d v$ continuous and $\eta^{2}$ reduced. Assume further
$(U_{j} \eta)^{2} \sim W_{j} \eta^{2}$, where $W_{j}$ are of the same character as $U_{j}$.
The condition 
$\frac{d U_{1} \eta}{d U_{2} \eta} \sim \frac{(U_{1} \eta)^{2}}{(U_{2} \eta)^{2}} \frac{d {}^{t} U_{1} v}{d {}^{t} U_{2} v} \neq const$, means that we must have $\frac{d {}^{t} U_{1} v}{d {}^{t} U_{2} v} \neq const$ except for a discrete set. 
Note that the argument depends on if $U_{j} \eta \rightarrow {}^{t} U_{j} v$ bijective.

 \section{Necessary condition on order in infinity}
 
 When $X dy - Y dx=0$, the direction $\frac{\xi}{\eta}$ associates the transformation 
$\xi \frac{\delta f}{\delta x} + \eta \frac{\delta f}{\delta y}$ to the point $(x,y)$ (\cite{Lie91}).
Assume $X(x,y) d y - Y(x,y) d x=0$ gives $\infty^{1}$ curves that represent an infinitesimal transform
$U f \equiv \xi(x,y) \frac{\delta f}{\delta x} + \eta(x,y) \frac{\delta f}{\delta y}$.  Given trajectories $\xi(x,y) = const$, we can using quadrature
determine an invariant set of trajectories (for $U f$) as $\eta(x,y)=const$. If we use $\xi,\eta$ instead
of $x,y$ in the differential equation, it gets a separated form and can be integrated using quadrature.

According to Lie, the class of $\infty^{1}$ curves $w(x,y)=const$
represent the infinitesimal transformation $Uf$, when $Uw$ is a function of $w$ alone, that is
$U w \equiv \xi \frac{\delta w}{\delta x} + \eta \frac{\delta w}{\delta y} = \Omega(w)$ (\cite{Lie91}).
Any projective transformation preserves separability (bijective and $const \rightarrow const$)

Congruent curves on $H_{m}$ are given by involution. In my model, we use $c \rightarrow c^{\bot_{E}}$,
which is not necessarily involutive. Assume congruence according to $c \rightarrow c^{\bot_{E}}$. 
Consider $\eta \rightarrow c$
and $v \rightarrow 1/c$, then when $c=0$, we have $x \bot_{E} y$, an oriented (one-sided) orthogonal.
The condition $1/ (\eta + i v) \rightarrow 0$ when either $x \bot_{E} y$ or $y \bot_{E} x$, defines a ``mean''
orthogonal. The domain for $(\eta,v)$ is such that $(\eta,v)$ has compact sub-level sets. When $(\eta,v)$ is semi-algebraic with compact sub-level sets, the set $\eta$
is semi-algebraic with cluster sets (or the set $v$)

\subsection{On collar points}

A linear dependence in $\infty$, for instance $d y / d x=const=\xi$, that reduces the order for 
the curve in $\infty$, means that the system is not integrable. Given a non-linear system,
there are at least two directions on $\Omega$ in a point (\cite{Lie96}). 
When $\frac{d y}{d x}-1 \equiv 0$ on a set of positive measure, we have $\frac{d}{d x}(y-x) \equiv 0$,
that is linear dependence in $\infty$. 

\newtheorem{lem9}[prop2]{Lemma}
\begin{lem9}
 Assume $\frac{d y}{d x}=\rho=e^{\phi}$ analytic with isolated singularities,
then $\{\zeta \quad \rho^{\bot_{E}}=0 \} \sim \{\zeta \quad \phi=0 \} \subset \Omega$ 
(geometric equivalence), where $\Omega$ pseudo-convex. 
When $\phi$ pseudo-convex, we can write $\Omega \sim \{ \rho^{\bot_{E}} > 0 \}$ locally. 
When $\phi \in L^{1}$, $\rho$ has algebraic singularities and $\rho^{\bot_{E}}$ has
isolated singularities.
\end{lem9}

Note that $(1 -\frac{d y}{d x}) \bigoplus \frac{d y}{d x}=1$ depends on the blow-up mapping.
The ideal $I(\rho^{\bot_{E}} =0)$ allows a global pseudo-base (using monotropy), but $I(\{ \rho = 0 \})$
has possibly only a local pseudo base. The first ideal is a geometric ideal over an algebraic set in $\zeta$.

The condition on collar points, is necessary for the congruence ($\bot_{E}$) to be well-defined. 
The condition for involution $\frac{d y'}{d x'}=\frac{d y}{d x}$, means
$\frac{d U_{1} y}{d U_{2} y} = \frac{d U_{1} x}{d U_{2} x}$.

Note that when $\gamma \in (I_{1})$ a rotation surface (symmetry), does not imply $U \gamma$
in a rotation surface, so we write $U \gamma= (x,y,z)$. When $\frac{d y}{d z}=\varphi(\frac{d x}{d z})$
and the branches are given by $\frac{d y}{d z}=const$, $\frac{d x}{d z}=const$. For the system $y d x - x d y=0$
the general form of $\Omega$ is $z =F(y/x)$, through bending of the main tangent to the region in question (\cite{Lie96}).

\subsection{Continuation of movement}
Consider the problem, when Schr\"odinger's model has an algebraic base. In this case we would give $U \gamma=
(x,y,z)$, where $z=P(x,y)$, for a polynomial $P$. Note also that $U \gamma$ algebraic does not imply $\gamma$ algebraic. Further $J U \gamma$
algebraic would be sufficient. When $F(U \gamma)=F(\gamma) + \widehat{V}$, according to Cousin's model
we do not have a normal model, but we may still have an algebraic base (\cite{Dahn13}).

Assume $y=e^{w(v)}$ and $x=e^{v}$. When $\frac{d}{d v} \log \eta \in L^{1}$ which implies $(\frac{d w}{d v})^{\bot_{E}} \in L^{1}$
and given an orthogonal division of $L^{1}$, $\frac{d w}{d v} \in L^{1}$. A sufficient condition for existence
of $w^{\bot_{E}}$ using annihilator theory, is that the range of $w$ is closed. 
 Assume $\log A=a$ and $\log AB \in L^{1}$.
When $a+ a^{\bot_{E}} \in L^{1}$, we have existence of $B$ such that $B=e^{a^{\bot_{E}}}$. 
Assume $\eta=e^{\phi}$ and $W \phi \in L^{1}$, then we have existence of $W^{\bot_{E}} \phi \in L^{1}$.
A continuation of $\eta$ can be written $\tilde{\eta}=e^{W^{\bot_{E}} \eta}$, through the collar point condition
we have $W \cap W^{\bot_{E}}=\{ 0 \}$. 

 \section{Dependence on movement}
\subsection{De-singularization}
 Assume $\gamma$ algebraic and $U$ a single movement, then $U \gamma$ is algebraic, but for composite 
movements, this is not necessarily true. For instance $P(z^{1/k})$ is not a polynomial in $z$. With the 
symbol representation $(F(x,y),M(\eta,1/\eta))$ and change of character of movement in the sense $M \sim {}^{t} M$, monodromy is not 
possible in $H$, but the movement can be determined in $\mathcal{D}'$. 

When $(I)=H_{m}$, for any $\gamma,\gamma'$ there is a $U$ such that $\gamma'=U \gamma$. Consider from
3-space $U_{j} \rightarrow U_{j}x=x(U_{j})$. When several 
movements are involved, $d x=\frac{d x}{d U_{1}}d U_{1} + \frac{d x}{d U_{2}} d U_{2} + \frac{d x}{d U_{3}} d U_{3}$
and $\frac{d x}{d U_{k}} / \frac{d x}{d U_{j}} \neq 0$, for $j \neq k$.  We assume $x$ is
defined on a cylindrical domain of $U_{j}$. When $x$ is polynomial, $\int_{\Omega_{U}} x d \mu=0$
implies $\mu(\Omega_{U})=0$ (measure zero). 

\newtheorem{def6}[prop2]{Definition}
\begin{def6}
 A regular approximation of a singular point, is generated by a movement $U$, such that the singular point is isolated on
the trajectory. If all rotations are regular, we could say that the point is isolated for rotation. If we
consider also the orthogonals, the limit could be regarded as ''two-sided``. The point is isolated
''globally`` if it is isolated for all movements.

\end{def6}

When $U^{\bot_{E}} \gamma \in H_{m}$, it is assumed analytic. When $U \gamma$ algebraic, we can assume that
$(U \gamma)^{\bot_{E}}$ analytic. 
Note that $\lambda_{1} < \mid \eta \mid < \lambda_{2}$ compact implies $\lambda_{1} < \mid \eta + i \frac{1}{\eta} \mid < \lambda_{2}$
compact. Assume $F(\gamma)(\zeta)=f(\zeta)$, where $f$ is reduced. In this case $f$ has compact sub-level sets.
Further, since $\mid F(x,\eta) \mid \leq C \mid \eta \mid$, this implies compact sub-level sets for $\eta$.
In the same manner, since $\mid F(y,v) \mid \leq C \mid v \mid$, it implies compact sub-level sets for $v$.
Thus there is no room for a spiral when $f$ is reduced. When $f$ is not reduced, we assume $F$ has non trivial kernel. 
When $f^{N}$ reduced, for some iteration index $N$, we are considering compact sub-level sets for the pair $(\eta,v)$
and we may have cluster sets for sub-level sets to $\eta$ in $\mbox{ker }F$ and $v$ in $\mbox{ ker }{}^{t} F$.

When $F$ corresponds to a Fredholm operator, we can consider $R(\gamma) \bigoplus Y_{0}$, where we assume $R(\gamma) \subset (I)$.
Thus $U \gamma =U_{0} \gamma + U_{1} \gamma$, where $F(U_{1} \gamma)=0$, which does not imply $U_{1} \gamma$ analytic.
When $F$ has a hypoelliptic representation over $\mbox{ ker }F$, it is sufficient to determine the continuation of movement,
that $U_{1} \gamma$ is continuous. More precisely, we assume $\{ U \gamma \} \subset (I)$ analytic with
$C$ discrete, and $\lim U \gamma \notin (I)$ continuous. 
 Example, $F(U_{S}^{\bot_{E}} \gamma)=\widehat{V} \sim 
\big[ F,U_{S}^{\bot_{E}} \big] (\gamma)$, with $\gamma$ algebraic. When $\eta$ has a limit in infinity,
the set $S=\{ \zeta \quad \eta + i v \leq \lambda \}$ is compact. 

Assume $\widehat{V} \in \mathcal{D}^{F '} (\Omega_{\zeta})$ and $\zeta_{t} \rightarrow 
\eta + i 1/\eta \rightarrow \widehat{V}$, where the last mapping is $\sim M U^{\bot_{E}}$. When we assume $U^{\bot_{E}} \eta + i 1/\eta$
has compact sub-level sets and when $J$ proper, mapping collar points on collar points, $V^{\bot_{E}} J \eta + i J 1/\eta$
has compact sub-level sets. When we assume $U_{1}^{\bot_{E}} \sim \tilde{U}_{2}$ where the movement changes character,
we can assume $F(\tilde{U}_{2} \gamma)=\widehat{V}$ and that $\gamma \rightarrow \tilde{U}_{2} \gamma$
preserves analyticity.
\subsection{Composite movements}
Consider as in the mixed model, $(I)=(I_{1}) \bigoplus (I_{1}^{\bot_{E}})$ and $\gamma'=A \gamma$ with
$\frac{d y_{1}'}{d x_{1}'}=\frac{d y_{1}}{d x_{1}}$,$\frac{d y_{2}'}{d x_{2}'}=\frac{d y_{2}}{d x_{2}}$.
In the case where $A$ is parabolic, we have $\frac{d y_{1}'}{d x_{1}'} - 1 \sim \frac{d y_{1}}{d x_{1}} - 1$.
Note that for the corresponding condition for $(I_{1}^{\bot_{E}})$, it is not sufficient to consider 
the tangents. 

Assume $\frac{d y}{d x}=\rho=Y/X$ analytic, that is $\{\zeta \quad \rho=0 \}$ analytic. If we assume $\rho \in L^{1}$
then also $\{ \zeta \quad \rho = 1 \}$ is analytic. When $\rho=e^{\varphi}$, $\{\zeta \quad \rho=1 \} \sim \{\zeta \quad \varphi=0 \}$, 
that is analytic when $\varphi$ is analytic.
Let $\rho^{\bot_{E}} \sim \rho - 1$ such that $\{ \rho^{\bot_{E}} = 0 \} \sim \{ \rho = 1 \}$.
In this case ``surjectivity means when $\rho \neq 0$, we have existence of limes of $\rho - 1$ as $x \rightarrow x_{0}$.
When $\rho \neq 0$, $\frac{d \rho}{d x}=0$ implies $\frac{d \varphi}{d x}=0$. We assume the ''transversals`` $\{\zeta \quad \varphi =0 \}^{\bot_{E}} \subset \{ \zeta \quad \rho \geq 1 \}$
semi-analytic. Assume $\psi=y-x$, then $\{ \zeta \quad \rho^{\bot_{E}} =0 \}$ $\simeq \{ \psi=\frac{d}{d x} \psi=0 \}$.
When $\psi=y - cx$, the set of singularities $y=cx,\rho=c$, correspond to $\rho$ on a scaled hyperboloid.

Assume $F(U \gamma)=F(\gamma) + \widehat{V}$ or $F(U^{\bot_{E}} \gamma)=\widehat{V}$.  
Assume $\{ P_{j} \}$ given points in $\gamma$ and $\zeta \in \Omega$,
then the involution condition is on the segments between the points. Assume the movement is $U=U_{1} + i U_{2}$. 
When $W$ defines a cylindrical domain in movement space and $V^{2} \sim \tilde{U}$, where $\tilde{U} \gamma \sim \mid U \gamma \mid^{2}$
we get the condition $2 \log V= \log \overline{U}^{\bot_{E}} + \log (-U)^{\bot_{E}}$. 

\newtheorem{def5}[prop2]{Definition}
\begin{def5}
 Assume $I$ an interval between $p,q$ and $I^{+}$ the interval for $U=U_{1}$ and $I^{-}$ the interval
for $U=U_{2}$. Assume $U$ the movement between $p,q$. Assume $U_{j}$ linearly
independent, that is $d U_{1}=d U_{2}$ and $U_{1}=U_{2}$ on $I$ implies $\mid I \mid=0$. Conversely,
$\mid I - I^{+} \mid=0$ implies $\mid (U - U_{1}) I \mid=0$ given that $(U - U_{1})$ a.c. Analogously,
$\mid (V - V_{1}) I \mid=0$ assuming $(V-V_{1})$ a.c. 

\end{def5}

Note that collar points for $U_{j}$ are mapped on
collar points for $U_{j}^{\bot_{E}}$, that is we have a proper mapping $\{ \mid \eta \mid < 1 \} \rightarrow 
\{ \mid v \mid < 1 \}$. Further,
collar points for $U_{j}$ is mapped by $J$ on to collar points for $V_{j}$. 

When $\eta^{2}$ reduced, $\eta^{2} d v \sim d \eta$ means that zero sets are preserved, when $\eta \rightarrow v$
bijective. When $\eta^{(N-2)/2}$ polynomial, the zero-sets are preserved. 
When $\eta^{N}$ reduced, consider $\eta'=\eta^{N/2}$ and $v'=v^{N/2}$, thus $\eta^{N} d v'=d \eta'$, 
then zero-sets are preserved when $\eta' \rightarrow v'$ bijective. 

\subsection{The transposed movement}
Assume $F(\gamma)(\zeta_{t})=F(U \gamma)(\zeta)=T F(\gamma)(\zeta)$. The scheme $\zeta_{t} \rightarrow T \rightarrow U$, is dependent
on if $T \rightarrow U$ is injective. 
We define $T F(\gamma)(\zeta)=F(\gamma)(\zeta_{t})$ as continuous on $C$ 
.When $\zeta_{t}$ 1-1 and closed, $T$ is surjective, thus given $F(\gamma)(\zeta_{t})$ on $C$, there
is a $\gamma_{1}$, such that $F(U \gamma_{1})=F(\gamma)(\zeta_{t})$. In order to come to a conclusion for 
$T \rightarrow U$, we must assume $F : H_{m} \rightarrow H_{m}$.  Assume $F$ preserves
dimension for the movement in $H_{m}$. When both $U$ and the orthogonal
movements are used to define singularities, we can conclude type of movement.

Consider $\Phi(\frac{d y}{d x})=\eta$ and $\frac{d \eta}{d x}> 0$, this means when $x$ real, that $\Phi(y)$ is convex.
When $\Phi(\frac{d y}{d x})={}^{t} \frac{d}{d x} \Phi(y)=\eta$ and $\Phi=const$, we have that
$(\Phi \frac{d}{d x} - \frac{d}{d x} \Phi)(y)=0$, for instance $\frac{d y}{d x}=\frac{y}{x}$, which
is the condition for determined tangent. Thus, when $\Phi=const$, we do not have a spiral approximation.

\newtheorem{lem11}[prop2]{Lemma}
\begin{lem11}
 Note that if we assume $\big[ I,U \big]=\big[ U,I \big]$, that is $U$ is algebraic,
then the action can be determined from phase space.

\end{lem11}

Assume $U^{\bot_{E}} \gamma=0$ implies $U^{\bot_{E}} e^{\phi_{n}} \sim e^{U^{\bot_{E}}_{1} \phi_{n}}=0$, 
then as $\mid \phi_{n} \mid \rightarrow \mid \phi \mid$, we have $U^{\bot_{E}} \gamma = g e^{U^{\bot_{E}}_{1} \phi}$
for $g$ measurable relative $U_{1}$ and $\log g + U_{1}^{\bot_{E}} \phi \in L^{1}$.

Consider $(\widehat{V},M(U \eta,v) - M(\eta,v))$, $({}^{t} \widehat{V},{}^{t} M(\eta,U v) - {}^{t} M(\eta,v))$ 
and finally $(\widehat{V}^{\bot_{E}},M(U^{\bot_{E}} \eta, v) - M(\eta,v))$. Thus when $M \sim {}^{t} M$, we must have $\widehat{V} \neq {}^{t} \widehat{V}$,
given that $U$ defines one simple movement. Further, when $\widehat{V}^{\bot_{E}}={}^{t} \widehat{V}$, we have an involutive
movement.
When $M$ is symmetric, its kernel is symmetric, that
is the reflection axes has a corresponding axes for invariance for $U^{\bot_{E}}$.

\subsection{Representation for change of movement}
On the hyperboloid, we can assume $U^{\bot_{L}}=I-U$ (orthogonal relative Lorentz metrics) and $U^{\bot_{L} \bot_{L}}=I -(I-U) = U$. 
 Note that $\frac{d F(U \eta)}{d U \eta} \frac{d U \eta}{d \eta}=0$
that is ''independent on $\eta$`` and $\frac{d U \eta}{d \eta} = 0$ implies $\frac{d F}{d U}=0$.
In the same manner for ${}^{t} F( U^{\bot_{E}} \gamma)={}^{t} \widehat{V}$ and $F(U^{\bot_{E} \bot_{E}} \gamma)=\widehat{V}^{\bot_{E}}$.
Note that $\frac{d U^{\bot_{E}} \eta}{d \eta}=1 - \frac{d U \eta}{d \eta}$, why if $\frac{d U^{\bot_{E}} \eta}{d \eta} < 0$ we have
$\frac{d U \eta}{d \eta} > 1 > 0$. 

When symmetry occurs in points,
the representation assumes that any movement in $\eta$ has a corresponding movement in $1/\eta$. When the representation
is reduced for translation in $\eta \rightarrow \infty$, it is reduced for rotation in $1/\eta \rightarrow 0$.  When $T F(\gamma)=F(U \gamma)$, we consider $T F - \lambda I$, where $\lambda$ is used
to adjust the reflection axes.

\newtheorem{def8}[prop2]{Definition}
\begin{def8}
 The condition $2 \log u_{3} \sim \log \overline{u}^{\bot_{E}} + \log (-u)^{\bot_{E}}$, where $u=(u_{1},u_{2})$
is taken in $\mathbf{R}^{3}$, that is we do not assume for instance $u_{1}^{\bot_{E} \bot_{E}}$ on the hyperboloid.
When $\tilde{u}_{1}^{\bot_{E}} \sim \tilde{u}_{2}$, we assume the planar objects (multivalent) from hyperbolic
geometry, corresponds to euclidean rotation. 

\end{def8}

Assume $\gamma=(x,y)$ a polynomial in $\zeta$, that is $y + i x(\zeta)$ and with the condition
$\mid \frac{x}{y} \mid \leq C$ for $\mid \zeta \mid \geq R$, where $C,R$ are constants. That is we assume
$\gamma$ is of real type. This condition in the hyperboloid, means that the movements are rotational.
Assume $V=\{ \zeta \quad x = c y \}$, then according to the real type condition $\zeta \in V$ implies $\mid \zeta \mid \leq R$,
that is $V \subset \subset \Omega$, where $\Omega$ is the domain for $\gamma$.

Note that when $C$ is given by constant surfaces to $F$, given $Y = - \frac{d F}{d x}$ Hamilton,
$C$ must be discrete, when $X,Y$ analytic. In the generalized moment problem, a continuous $C$
implies $F \in \mathcal{D}' \backslash H$.

\subsection{Dependence of measure for movement}
Assume associated to $(U_{1},U_{2},U_{3})$ in a cylindrical domain, a measure $\mu$ on the web $\Omega$ of the cylinder.
That is $U_{1}^{2} + U_{2}^{2}=U_{3}^{2}$, and the measure on simple movements zero, $\mu \sim \mu_{1} \otimes \mu_{2}$.
Note that $\mu(U_{1},0,0) \cap \Omega$ and $\mu(0,U_{2},0) \cap \Omega$ are points.  
Using Lie's involution condition $\frac{d U y}{d U x}=\frac{d y}{d x}=\rho$ regular. 

When $\eta \in \mathcal{D}_{L^{1}}' (\mathcal{U})$, we assume $U \eta(\zeta) = \eta(\zeta_{t})$, where
$t$ is parameter for the movement $U$. 
When $\eta$ is analytic in $\zeta_{t}$,
we have that $\eta$ is analytic in $t$. Note that in this case, singularities are locally algebraic, 
when $\log \mid \eta (\zeta_{t}) \mid \in L^{1}$.

\newtheorem{lem12}[prop2]{Lemma}
\begin{lem12}
 According to the above,
when $M_{1}(\gamma)=\int \gamma dU_{1}$, we have $M_{2}(\gamma)=\int \gamma/\alpha d U_{2}$ and $M_{3}(\gamma)=\int \gamma/\beta d U_{3}$.
When $U=U_{1} + i U_{2}$ (on a cylindroid), we consider $\tilde{U_{j}} \sim (U_{j})(-U_{j})$,
$j=1,2$ and note that $\tilde{U_{1}} \sim \tilde{(-U_{1})}$
and $(\tilde{U_{1}} + i \tilde{U_{1}}^{\bot_{E}})^{\bot_{E}} \sim -i (\tilde{U_{1}} + i \tilde{U_{2}})$, where we used
that $\tilde{U_{2}}^{\bot_{E}} \sim \tilde{U_{1}}$. On the set $C$ when $\frac{d U_{1}}{d U_{2}} \equiv 1$ we have $\gamma=\gamma_{0}$.
This means that the cylindroid movement preserves analyticity
on the set $C$ where the movement changes character.
When $M \sim {}^{t} M$ in euclidean scalar product, we can define $\gamma$ such
that $\gamma/\alpha \bot_{E} H$ ($H$ being holomorphic functions). 
\end{lem12}
Assume $\Omega_{j}$ a local neighbourhood generated by a movement $U_{j}(t)$ and $\frac{d U_{2}}{d U_{3}}=\beta$
and $\frac{d U_{2}}{d U_{1}}=\alpha$, then $M_{2}(\gamma)=\int_{\Omega_{2}} \gamma \frac{1}{\alpha} d U_{2}$
and $M_{3}(\gamma)=
\int_{\Omega_{3}} \gamma \frac{\beta}{\alpha} d U_{3}$.

Concerning reflection, $(U_{1},U_{2}) \rightarrow (U_{1},-U_{2})$ is planar reflection. 
Note that $\tilde{U}_{2}=(-U_{2}) U_{2}$, that is $\tilde{U}_{2}^{\bot_{E}}=-\tilde{U}_{1}$ and $(\tilde{U}_{1} + i \tilde{U}_{2})^{\bot_{E}}=
(\tilde{U}_{1} - i \tilde{U}_{1})^{\bot_{E}} = (\tilde{U}_{2} - i \tilde{U}_{1})$, that is the mapping over $C$
is pure.

Assume $M_{C}(\frac{1}{y})=\int_{C} \frac{dy}{y}=\int_{C} \rho_{n} \frac{d x}{y}$, where $\rho_{n}=\frac{d y}{d x} \mid_{I_{n}}$
and $\Phi(\rho_{n})=\eta_{n}$ dependent on division element $I_{n}$. Thus, when $\rho_{n}$ analytic,
it is locally 1-1, as $n \rightarrow \infty$. Note when $L=\{ (x,y) \quad \eta_{n}=c \neq 0 \}$,
$L'=\{ (x,y) \quad \rho_{n}= c' \neq 0 \}$ and $L''=\{ (x,y) \quad \frac{d^{2} y}{d x^{2}}=c'' \neq 0 \}$.
Thus, when $0 \in L$, $\eta=\frac{0}{0}$ and $0 \notin L' \cap L''$ (cf. very regular boundary).
When $\rho_{n}=\frac{d y}{d x}$ dependent on division, the index is not dependent on division.
We have $M_{C_{n}} \mid I_{n} \mid u_{1} + M_{C_{n}} \mid I_{n} \mid u_{2} \rightarrow 0$, as $\mid I_{n} \mid \rightarrow 0$.

\subsection{Dependence on spectrum for movement}
Consider $<U \gamma,{}^{t} F \psi >=<\gamma, \lambda \psi>$, then when $U \gamma=\gamma$,
we have ${}^{t} F \psi-\lambda \psi \bot_{E} \gamma$. If we continue this equation to $\tilde{F}$ localizer
to a reduced operator, that is with trivial kernel and surjective, where ${}^{t} F \psi = \lambda \psi$
we have $ U \gamma= \gamma$., otherwise we consider $\psi \notin \mbox{ ker } {}^{t} F$ and when ${}^{t} F \psi = \lambda \psi$,
$\psi \bot_{E} U \gamma - \gamma$.

\newtheorem{lem13}[prop2]{Lemma}
\begin{lem13}
 Let $T F(\eta)=F(U^{\bot_{E}} \eta)=\lambda F(\eta)$, then when $\lambda \eta \bot_{E} F$ iff $\mid \eta \mid < 1$
($U = U_{1}$) or $\eta \bot_{E} F$ iff $\mid \eta / \lambda \mid < 1$ iff $\mid \eta \mid < \lambda$, we can determine
the movement.
In particular, when $\ker F = \{ 0 \}$ (bijective) and $U^{\bot_{E}} \eta =0$ implies $\mid \eta \mid < \lambda$,
the movement is determined. Further we can regard $F(\eta)$ as an eigenvector to $T - \lambda I$.

\end{lem13}

Given $U$ singular, there is a line through $0$, such that $U^{\bot_{E}} \eta =0$ implies $\mid \lambda \eta \mid < 1$, for instance
$U = U_{1}$. If $1 \in \sigma_{p}(U)$ that is $U^{\bot_{E}} \eta=0$, then $\eta$ is a vector. If $1 \in \sigma_{c}(U)$
that is $U^{\bot_{E}} \eta =0$, then $\eta$ is a point. Alternatively, consider $T F(\gamma)=T_{p} F(\gamma)$, where
$T_{p}$ is a parabolic movement. We have that $L(U_{p} x, U_{p} x)=L(x + \chi l,x + \chi l)$, where $L(l,l)=0$
corresponds to scaling of the hyperboloid, $L(x + \chi l,x + \chi l)=\lambda L(x,x)$ (\cite{Riesz43}). When $T F(\eta)= \lambda T_{p} F(\eta) = \lambda F(\eta)$,
implies $\mid \lambda \eta \mid = 1$ and as $\lambda \rightarrow \infty$, this implies
$\mid \eta \mid < 1$, which implies $T=U_{1}$. 

\subsection{Essential spectrum}
We consider $\widehat{V} \in \mbox{ ker }E$ and we define $v_{ess}$ as 
the set of parameters where $E(U^{\bot_{E}} \gamma)$ is not closed. 
Assume $E - I \in \widehat{C}^{\infty}$ on a Banach-space,
that is $E$ is assumed corresponding to very regular action. Then $E$ has closed range and the dimension for $N(E)$ is finite.
Thus in the case when $\mbox{ ker }E$ non-trivial, we have $X=R(E) \bigoplus Y_{0}$
and $E$ is not considered as closed on $Y_{0}$. This motivates $v_{ess}$ as a correspondent to $\sigma_{ess}$.
The $\sigma_{ess}$ is spectrum for the operator modulo regularizing action. 

\newtheorem{def1}{Definition}[section]
\begin{def1}
 Define an analogue to $\sigma_{ess}$ according to $v_{ess}=$
$\{ \lambda$ such that $\mbox{ ker }E_{\lambda}  \neq \{ 0 \} \}$ 
where $E_{\lambda}$ corresponds to $\gamma - \lambda$ (or $P-\lambda$ when $PE(f) \sim f$)

\end{def1}

Note that given an operator $A$ Fredholm, then for $\lambda \in \sigma_{ess}(A)$, we have that
$\mbox{ index }(A)$ is constant on compact sets. When the movement is spiral, we have sub-level sets with
cluster-sets in both ends. The limit for $(\eta,v)$ when $v \sim 1/\eta$, is symmetric in the sense that $(\eta,v)$ has compact sub-level sets
iff $(v,\eta)$ has compact sub-level sets. However $I \prec \prec \eta + i v$ does not imply $I \prec \prec \eta$
or $I \prec \prec v$.

When $\widehat{V}=0$ implies $\mid \eta \mid \leq 1$, we have one movement, otherwise we divide the domain
into segments of ``constant sign''. When $F(U \gamma) - F(\gamma)=\widehat{V}$
in $L^{1}$, we can use Radon-Nikodym to represent the movement by a function in $L^{1}$, that is 
$\widehat{V}= F((g-1) \gamma)$ $g \in L^{1}$ ($F$ linear in $\gamma$).
Note that $U^{\bot_{E}} \gamma$ analytic, gives a continuous mapping to $\zeta_{t}$, that is $\{ \mid \eta \mid < 1 \}$
has a correspondent set $\Omega_{\eta}=\{ \zeta_{t} \}$ and in the same manner for $\Omega_{v}$. Assume $\widehat{V}=0$ implies $\mid \eta \mid < 1$ and
$\Omega_{\eta}$ non compact and $\Omega_{v}$ compact. When ${}^{t} \widehat{V}=0$ implies $\mid v \mid < 1$
implies $\Omega_{\eta}$ compact and $\Omega_{v}$ non compact.

\subsection{Modulo parabolic movements}
Let $G(\eta)=F(x,\eta)$ and denote $U \eta$ for $U \gamma \rightarrow \frac{U y}{U x}$.
Assume $T G(\eta)=G(U^{\bot_{E}} \eta)$, then modulo parabolic movements, when $T G(\eta)=\lambda G(\eta)$, implies
$\mid \lambda \eta \mid < 1  $, we have for large $\lambda$, the movement must be translation. Symmetric points
for $M(\eta,1/\eta)$ are not possible. For a single movement $M \neq {}^{t} M$. When $\lambda \sim 1$, we may have change of character for the movement.
 Note that when the movement is monotonous, it is not necessarily single.

\newtheorem{lem14}[prop2]{Lemma}
\begin{lem14}
 Concerning $\big[ U_{3},U_{1} \big]$, we claim that this movement has the same character as $U_{1}$.

\end{lem14}

When $T G(\gamma) = \lambda T_{p} G(\gamma)$ and we consider $(F, M(\eta,1/\eta))$ with $M \sim {}^{t} M$,
then for $\mid \eta \mid=1$ $M(\eta, \lambda 1/\eta)=\overline{\lambda} M(\eta,1/\eta)$, that is when
${}^{t} M \sim M$, $\lambda$ must be real. 

The continuous composite
movement, can be approximated infinitesimally by sequential movements. That is $\Sigma I_{j} - I \rightarrow 0$
for parameter intervals, implies $\Sigma ( U_{j} - U) \eta \rightarrow 0$, assuming the movement monotonous.
The limit is not necessarily a reflection, that is there is not an axes to the movement. When $U_{1}=\big[ U_{1},U_{3} \big]$
we can consider $G(U_{1}^{\bot_{E}} \eta)=\lambda G(\eta)$ where $U_{3} \eta \sim \lambda \eta$, why
$\mid \lambda \eta \mid < 1$, that is for large $\lambda$ the movement is a single translation.
Note that the micro-local analysis is completely determined by translational movements. 

When the movement does not change character as $\lambda \rightarrow \infty$, the movement is simple.
Assume existence of a point $p \in C$, such that $\mid \lambda \eta \mid < 1$ when $\lambda < \lambda_{1}$
and $\mid \lambda \eta \mid > 1$ as $\lambda_{1} < \lambda$, then we see that $\sigma$ can be used 
to determine the eigenvectors corresponding to change of character.

Consider the example
with joint spectrum, $\lambda M(\eta,v)=\lambda_{1} \lambda_{2} M(\eta,v)=M(\lambda_{1} \eta, \lambda_{2} v)$.
When $\lambda < 1$, we can chose $\lambda_{1}$ large 
such that when $\mid \lambda_{1} \eta \mid < 1$ and $\mid \lambda_{2} v \mid < 1$, we have 
$\lambda_{2} < \mid \eta \mid < 1/\lambda_{1}$
and when $1 < \lambda$ such that $\mid \lambda_{1} v \mid < 1$ and $\mid \lambda v \mid < 1$, 
we can assume $\lambda_{2} < \mid \eta \mid < 1/\lambda_{1}$. Change of character is possible
when $\lambda \sim 1$.

\bibliographystyle{amsplain}
\bibliography{schrodinger}

\end{document}